\documentclass[12pt]{amsart}
\usepackage{amsmath,amsthm}
\usepackage{amssymb}
\usepackage{euscript}
\DeclareMathAlphabet{\EuRm}{U}{eur}{m}{n}
\SetMathAlphabet{\EuRm}{bold}{U}{eur}{b}{n}
\begin{document}
%
%
\newcounter{numb}
\newcounter{num}
%
%
\swapnumbers
\newtheorem{thm}{Theorem}[section]
\newtheorem*{tha}{Theorem A}
\newtheorem{lemma}[thm]{Lemma}
\newtheorem{prop}[thm]{Proposition}
\newtheorem{cor}[thm]{Corollary}
\theoremstyle{definition}
\newtheorem{defn}[thm]{Definition}
\newtheorem{example}[thm]{Example}
\newtheorem{notation}[thm]{Notation}
\newtheorem{summary}[thm]{Summary}
\newtheorem{fact}[thm]{Fact}
\theoremstyle{remark}
\newtheorem{remark}[thm]{Remark}
\newtheorem{assume}[thm]{Assumptions}
\newtheorem{note}[thm]{Note}
\newtheorem{ack}[thm]{Acknowledgements}

\numberwithin{equation}{section}
%
%
\def\sect{\setcounter{thm}{0} \section}
%
%
\newcommand{\xra}[1]{\xrightarrow{#1}}
\newcommand{\lra}{\longrightarrow}
\newcommand{\epi}{\to\hspace{-5.2mm}\to}
\newcommand{\hra}{\hookrightarrow}
\newcommand{\rlh}{\rightleftharpoons}
\newcommand{\oulra}[2]{\substack{{#1}\\ \rlh \\ {#2}}}
\newcommand{\lle}{<\hspace{-1.5mm}<}
\newcommand{\hsp}{\hspace{10 mm}}
\newcommand{\hs}{\hspace{5 mm}}
\newcommand{\hsm}{\hspace{2 mm}}
\newcommand{\vs}{\vspace{7 mm}}
\newcommand{\vsm}{\vspace{2 mm}}
\newcommand{\rest}[1]{\lvert_{#1}}
\newcommand{\DEF}{:=}
%
%
\newcommand{\lrb}[2]{[\![{#1},{#2}]\!]}
\newcommand{\lrbb}[1]{[\![{#1}]\!]}
\newcommand{\lrc}[1]{\langle{#1}\rangle}
\newcommand{\llrr}[1]{\langle\!\langle{#1}\rangle\!\rangle}
\newcommand{\lro}[1]{[{#1}]}
\newcommand{\lrp}[1]{({#1})}
\newcommand{\lrv}[1]{\lvert{#1}\rvert}
\newcommand{\lamt}[3]{\lambda_{3}({#1}\otimes{#2}\otimes{#3})}
\newcommand{\lamf}[4]{\lambda_{4}({#1}\otimes{#2}\otimes{#3}\otimes{#4})}
\newcommand{\eps}[2]{(-1)^{\lrv{#1}\lrv{#2}}}
\newcommand{\epsi}[2]{(-1)^{\lrv{#1}\lrv{#2}+1}}
\newcommand{\npr}[2]{[{#1},{#2}]}
%
%
\newcommand{\colim}{\operatorname{colim}}
\newcommand{\diag}{\operatorname{diag}}
\newcommand{\db}{\operatorname{db}}
\newcommand{\gr}{\operatorname{gr}}
\newcommand{\ho}{\operatorname{ho}}
\newcommand{\Hom}{\operatorname{Hom}}
\newcommand{\Image}{\operatorname{Im}\,}
\newcommand{\Ker}{\operatorname{Ker}}
%
%
\newcommand{\A}{{\mathcal A}}
\newcommand{\Alg}{{\EuScript Alg}}
\newcommand{\C}{{\mathcal C}}
\newcommand{\Do}{\mathbf{\Delta}^{op}}
\newcommand{\FF}{{\mathcal F}}
\newcommand{\II}{\EuScript{I}}
\newcommand{\JJ}{{\mathcal J}}
\newcommand{\HH}[1]{\EuScript{H}^{(#1)}}
\newcommand{\KK}{\EuScript{K}}
\newcommand{\LL}{{\mathcal L}}
\newcommand{\Lie}{{\EuScript Lie}}
\newcommand{\Lc}{\LL_{0}}
\newcommand{\M}{{\EuScript M}}
\newcommand{\Pa}{$\Pi$-algebra}
\newcommand{\Set}{{\EuScript Set}}
\newcommand{\Ss}{{\mathcal S}}
\newcommand{\Sa}{\Ss_{\ast}}
\newcommand{\TT}{{\mathcal T}}
\newcommand{\TTs}{\TT_{\ast}}
\newcommand{\TTo}{\TT_{1}}
\newcommand{\TTq}{\TT_{\Q}}
\newcommand{\VV}{{\mathcal V}}
\newcommand{\Vect}{{\EuScript Vec}}
%
%
\newcommand{\bA}[1]{{\mathbb A}\lrc{#1}}
\newcommand{\bJ}[1]{{\mathbb J}\lrc{#1}}
\newcommand{\bL}[1]{{\mathbb L}\lrc{#1}}
\newcommand{\N}{\mathbb N}
\newcommand{\Q}{\mathbb Q}
\newcommand{\R}{\mathbb R}
\newcommand{\bV}[1]{{\mathbb V}\lrc{#1}}
\newcommand{\ZZ}{\mathbb Z}
%
%
\newcommand{\bface}[1]{\mathbf{#1}}
\newcommand{\LLL}{\bface{L}}
\newcommand{\bS}[1]{\bface{S}^{#1}}
\newcommand{\X}{\bface{X}}
\newcommand{\XQ}{\X_{\Q}}
%
%
\newcommand{\Hi}[1]{H'_{#1}}
\newcommand{\His}{\Hi{\ast}}
\newcommand{\pis}{\pi_{\ast}}
\newcommand{\pii}[1]{\pi^{i}_{#1}}
\newcommand{\pie}[1]{\pi^{e}_{#1}}
\newcommand{\pim}{\pi_{\ast-1}}
\newcommand{\piq}[1]{\pi_{\ast}{#1}_{\Q}}
%
%
\newcommand{\As}{A_{\ast}}
\newcommand{\Ass}{A_{\ast,\ast}}
\newcommand{\hAss}{\hat{A}_{\ast,\ast}}
\newcommand{\Bs}{B_{\ast}}
\newcommand{\Bss}{B_{\ast,\ast}}
\newcommand{\Cs}{C_{\ast}}
\newcommand{\Css}{C_{\ast,\ast}}
\newcommand{\Hss}{H_{\ast,\ast}}
\newcommand{\Js}{J_{\ast}}
\newcommand{\Jss}{J_{\ast,\ast}}
\newcommand{\Ks}{K_{\ast}}
\newcommand{\Ls}{L_{\ast}}
\newcommand{\LX}{L_{X}}
\newcommand{\Lx}[1]{\Pi_{#1}^{X}}
\newcommand{\Lsx}{\Lx{\ast}}
\newcommand{\hLx}{\hat{L}_{X}}
\newcommand{\Lss}{L_{\ast,\ast}}
\newcommand{\hLss}{\hat{L}_{\ast,\ast}}
\newcommand{\hLs}{\hat{L}_{\ast}}
\newcommand{\Ms}{M_{\ast}}
\newcommand{\Ts}{T_{\ast}}
\newcommand{\Vs}{V_{\ast}}
\newcommand{\Xs}{X_{\ast}}
\newcommand{\Xss}{X_{\ast\ast}}
\newcommand{\Ys}{Y_{\ast}}
%
%
\newcommand{\Ad}{A_{\bullet}}
\newcommand{\Add}{A_{\bullet\bullet}}
\newcommand{\Ads}{A_{\bullet,\ast}}
\newcommand{\Bd}{B_{\bullet}}
\newcommand{\Cd}{C_{\bullet}}
\newcommand{\Cdd}{C_{\bullet\bullet}}
\newcommand{\Cds}{C_{\bullet,\ast}}
\newcommand{\Cu}[3]{C_{{#2},{#3}}^{(#1)}}
\newcommand{\Cus}[2]{\Cu{#1}{#2}{\ast}}
\newcommand{\Cuds}[1]{\Cus{#1}{\bullet}}
\newcommand{\Cnd}{C_{n,\bullet}}
\newcommand{\CdD}{\Cdd^{\Delta}}
\newcommand{\hCd}{\hat{C}_{\bullet\bullet}^{\Delta}}
\newcommand{\Eds}{E_{\bullet,\ast}}
\newcommand{\Jd}{J_{\bullet}}
\newcommand{\Jdd}{J_{\bullet\bullet}}
\newcommand{\Jds}{J_{\bullet,\ast}}
\newcommand{\Ju}[3]{J_{{#2},{#3}}^{(#1)}}
\newcommand{\Jus}[2]{\Ju{#1}{#2}{\ast}}
\newcommand{\Juds}[1]{\Jus{#1}{\bullet}}
\newcommand{\Kd}{K_{\bullet}}
\newcommand{\Ld}{L_{\bullet}}
\newcommand{\Lds}{L_{\bullet,\ast}}
\newcommand{\Uds}{U_{\bullet,\ast}}
\newcommand{\Vd}{V_{\bullet}}
\newcommand{\Vds}{V_{\bullet,\ast}}
\newcommand{\Vss}{V_{\ast,\ast}}
\newcommand{\Wd}{W_{\bullet}}
\newcommand{\Wds}{W_{\bullet,\ast}}
\newcommand{\Wdd}{W_{\bullet\bullet}}
\newcommand{\Xdd}{X_{\bullet\bullet}}
\newcommand{\Xd}{X_{\bullet}}
\newcommand{\Yd}{Y_{\bullet}}
\newcommand{\Ydd}{Y_{\bullet\bullet}}
\newcommand{\Zd}{Z_{\bullet}}
%
%
\newcommand{\F}[2]{{\EuScript F}^{#1}(#2)}
%
%
\newcommand{\DD}[1]{{\EuScript D}}
\newcommand{\Di}[2]{\DD^{#1}_{(#2)}}
\newcommand{\Sph}[2]{{\EuScript S}^{#1}_{(#2)}}
%
%
%
\title{Homotopy operations and rational homotopy type}
\author{David Blanc}
\address{Dept.\ of Mathematics, University of Haifa,
Haifa 31905, Israel}
\email{blanc@math.haifa.ac.il}
\date{October 20, 2001; revised August 31, 2002}
\maketitle
%
%
\sect{Introduction}
\label{ci}

In \cite{HStaO} and \cite{FelDT} Halperin, Stasheff, and F\'{e}lix showed how 
an inductive\-ly-defined sequence of elements in the cohomology of a graded
commutative algebra over the rationals can be used to distinguish among the 
homotopy types of all possible realizations, thus providing a collection of
algebraic invariants for distinguishing among rational homotopy types of spaces.
There is also a dual version, in the setting of graded Lie algebras 
(see \cite{OukiH}).

However, these authors provided no homotopy-theoretic interpretation of these 
invariants, which are defined in terms of differential graded algebras 
(resp.\ Lie algebras) and their possible perturbations. 

The goal of this paper is to provide such an interpretation, 
in terms of higher rational homotopy operations, and thus to make
sense of the following

\begin{tha}
For any simply-connected space $\X$, there is a sequence of higher homotopy
operations taking value in \ $\pis\X$, \ which, together
with the rational homotopy Lie algebra \ $\pim\X\otimes_{\ZZ}\Q$ \ itself, 
determine the rational homotopy type of $\X$. 
\end{tha}
\noindent (See Theorem \ref{tone} below). At the same time, we provide
a more concrete (rational) version of the general theory of higher
homotopy operations provided in \cite{BMarkH}.

It should be noted that an integral version of the Lie algebra case is contained
in \cite{BlaAI} (see also \cite{BGoeR,BDGoeC}), and the mod $p$ homology analogue of the 
Halperin-Stasheff-F\'{e}lix theory appears in \cite{BlaS}. Moreover, in
\cite{BlaHH} we showed that the (integral) homotopy type of a space $\X$ is 
in fact determined by its homotopy groups \ $\pis\X$, \ together with the action 
of all primary homotopy operations on it, and of certain higher homotopy 
operations (see \cite{BlaHO,BlaC} for subsequent modifications and improvements).

However, if we are interested only in the \textit{rational} homotopy type of a 
simply-connected space $\X$, Whitehead products are the only
non-trivial primary homotopy operations on the rational homotopy 
groups \ $\piq{\X}=\pis\X\otimes\Q$, \ which, after re-indexing, constitute a 
graded Lie algebra over $\Q$. \ The relevant higher order
operations are also simpler than in the integral case.
Thus we hope that the rational version of this theory will be both easier to 
understand, and more accessible to computation.

Moreover, the higher operations we define are certain subsets of \ 
$\pis\X$, \ indexed by elements in homology groups of a certain inductively 
defined collection of differential graded Lie algebras (DGLs) defined below, so
we provide an explicit correspondence between our higher operations and the
corresponding elements in cohomology groups with coefficients in \ $\pis\X$
(provided by the Halperin-Stasheff-F\'{e}lix theory) \ -- \ a correspondence 
which was lacking in the integral case.

Finally, while a notion of higher homotopy operations for a differential 
graded Lie algebra $L$ has been defined in the special case of higher Whitehead 
products (also known as ``Lie-Massey products'' \ -- \ see 
\cite{AlldR,AlldR2,AArkS,RetaMO,RetaL,TanrH}), in general it is not clear how to 
represent such rational operations  as \textit{integral} higher 
order operations in \ $\pis \X$, \ if $L$ represents the rational homotopy type 
of a topological space $\X$. In order to address this problem, we must consider 
a somewhat ``flabbier'' model of rational homotopy than that provided 
by differential graded Lie algebras, namely a certain class of 
differential graded non-associative algebras (see section \ref{cna} below).

Thus we also provide a (somewhat incomplete) answer to the following
question: what additional structure on the ordinary homotopy 
groups \ $\pis\X$ \ of a simply-connected space $\X$, \ beyond the
Whitehead products, is needed to determine its homotopy type up to
rational equivalence?

\subsection{Notation and conventions}
\label{snac}\stepcounter{thm}

The ground field for all vector spaces, algebras, and tensor products
will be $\Q$ (the rationals), unless otherwise stated.

$\TTs$ \ denotes the category of pointed $CW$ complexes
with base-point preserving maps, and by a \textit{space} we shall always
mean an object in $\TTs$, \ which will be denoted by a boldface
letter: \ $\X, \bS{n},\dotsc$.\  The subcategory of $1$-connected spaces
is denoted by \ $\TTo$, \ and the rationalization of a space \ $\X\in\TTo$ \ 
is \ $\XQ$. \ The category of rational $1$-connected topological 
spaces is denoted by \ $\TTq$.

Let $\Delta$ denote the category of ordered sequences \
$\mathbf{n}= \langle 0,1,\dotsc,n \rangle$ \ ($n\in \N$), \ with 
order-preserving maps.\hsm 
For any category $\C$, we let \ $s\C$ \ denote the category of
simplicial objects over $\C$  \ -- \ i.e., functors \ $\Do\to\C$ \ 
(cf.\ \cite[\S 2]{MayS}); objects therein will be written \
$\Ad,\dotsc$. \ If we omit the degeneracies, we have a 
$\Delta$-\textit{simplicial} object, which we denote by \ 
$\Ad^{\Delta},\dotsc$.

The category of non-negatively graded objects over a category 
$\C$ will be denoted by \ $\gr\C$, \ with objects written \ $\Ts,\dotsc$; \ 
we will write \ $\lrv{x}=p$ \ if \ $x\in T_{p}$. \ 
An upward shift by one in the indexing will be denoted by \ 
$\Sigma:\gr\C\to\gr\C$, \ so that \ $(\Sigma\Xs)_{k+1}=X_{k}$, \ and \ 
$(\Sigma\Xs)_{0}=0$. \ The category of graded vector spaces is denoted
by $\VV$.

The category of chain complexes (over $\Q$) will be denoted by \
$d\VV$, \ and that of double chain complexes by \ $dd\VV$. \ 
The differential of any differential graded object is written $\partial$ \ (to
distinguish it from the face maps \ $d_{i}$ \ of a simplicial object).

If $\C$ is a closed model category \ (cf.\ \cite[I]{QuH} or 
\cite[II, \S 1]{QuR}), \ we denote by \ $\ho\C$ \ the corresponding
homotopy category. \ If \ $X\in\C$ \ is cofibrant and \ $Y\in\C$ \ is
fibrant, we denote by \ $[X,Y]_{\C}$ \ the set of homotopy classes of
maps between them.

Let \ $\Set$ \ denote the category of sets, \ $\Vect$ \ the category of
vector spaces (over $\Q$), \ $\Lie$ \ the category of Lie algebras, 
and \ $\Alg$ \ the category of non-associative algebras. \ 
We write \ $\Ss$ \ rather than \ $s\Set$ \ for the category of
simplicial sets, \ and \ $\Sa$ \ for the category of \textit{pointed}
simplicial sets.

\subsection{Organization:}
\label{sorg}\stepcounter{thm}

In section \ref{clm} we review some background material on the Quillen 
DGL model for rational homotopy theory, and describe a bigraded variant of it; 
and in section \ref{csr} we give some more background on simplicial 
resolutions. 
These are applied to the rational context in section \ref{crrs},
where we also define higher order homotopy operations for DGLs. These
appear as the obstructions to realizing certain algebraic
equivalences, and serve to determine the rational homotopy type of a
simply-connected space. We give a first approximation to Theorem A in
\S \ref{shho}.

In section \ref{cmr} we explain how to translate the usual bigraded
and filtered DGL models into simplicial DGLs, which allows us to
construct appropriate minimal simplicial resolutions.
In section \ref{chd} we define the homology and cohomology of a DGL
(after Quillen), and show that the obstructions we define above 
actually take value in the appropriate cohomology groups. 
Finally, in section \ref{cna} we describe a non-associative
differential graded algebra model for rational homotopy theory, which 
facilitates the translation of the higher homotopy operations
described above into integral homotopy operations.
We summarize our main results in Theorem \ref{tone}.

\begin{ack}\label{aa}\stepcounter{subsection}
I would like to thank Ron Livn\'{e} for arousing my interest in
rational higher homotopy operations, Martin Arkowitz and
Jean-Michel Lemaire for providing me with copies of hard-to-get
theses, and Steve Halperin for some useful comments. I would also like
to thank the referee for his remarks and corrections.
\end{ack}

%
%
\sect{Lie models}
\label{clm}

In this section we briefly recall some well-known definitions and
facts of rational homotopy theory, and describe variants thereof.

\subsection{Differential graded Lie algebras}
\label{sdgl}\stepcounter{thm}

Let $\LL$ \ denote the category of graded Lie algebras, or \textit{GL}'s. \ 
An object \ $\Ls\in\LL$ \ is thus a graded vector space: \ 
$\Ls=\oplus_{n=0}^{\infty} L_{n}$ \ over $\Q$, equipped with a bilinear graded 
product \ $[~,~]:L_{p}\otimes L_{q}\to L_{p+q}$ \ for each \ $p,q,\geq 0$, \ 
such that \ $[x,y] = \epsi{x}{y}[y,x]$ \ and \ 
$\eps{x}{z}[[x,y],z] + \eps{y}{x}[[y,z],x] + \eps{z}{y}[[z,x],y]=0$. \ 
We denote by \ $\Lc$ \ the full subcategory of all \emph{connected} 
graded Lie algebras \ -- \ that is, those for which \ $L_{0}=0$.

The free graded Lie algebra generated by a graded set \ $\Xs$ \ is
denoted by \ $\bL{\Xs}$. \ The functor \ ${\mathbb L}:\gr\Set\to\LL$ \
is left adjoint to the forgetful ``underlying graded set'' functor \ 
$U:\LL\to\gr\Set$, \ and it factors through $\VV$: \ that is, \ 
$\bL{\Xs}=L(\bV{\Xs})$, \ where \ $\bV{\Xs}\in\VV$ \ is the graded 
vector space with basis \ $\Xs$ \ and \ $L(\Vs)$ \ is the free Lie
algebra on the graded vector space \ $\Vs$ \ (defined as the
appropriate quotient of the graded tensor algebra).

The category of \textit{differential graded Lie algebras}, or \textit{DGL}s, 
will be denoted by \ $d\LL$, \ with \ 
$d\Lc$ \ the subcategory of \textit{connected} Lie algebras (i.e., those
with \ $L_{0}=0$). \ An object \ $L=\lrp{\Ls,\partial_{L}}\in d\LL$ \ is a
graded Lie algebra \ $\Ls\in\LL$, \ together with a differential \ 
$\partial_{L}=\partial^{n}_{L}:L_{n}\to L_{n-1}$, \ for each \ $n>0$, \ 
such that \ $\partial^{n-1}_{L}\circ \partial^{n}_{L}=\{0\}$ \ and \ 
$\partial_{L}[x,y]= [\partial_{L}x,y] + (-1)^{\lrv{x}}[x,\partial_{L}y]$.

The homology of the underlying chain complex of a DGL \
$L=\lrp{\Ls,\partial}$ \ will be denoted \ $\His L$, \ to distinguish it
from the DGL homology defined in \S \ref{dhdgl} below. Because the 
differential $\partial$ is a derivation, \ $\His L$ \ inherits from $L$ the 
structure of a graded Lie algebra. 

A morphism of DGLs which induces an isomorphism
in homology will be called a \textit{quasi-isomorphism}, or \textit{weak 
equivalence}, \ denoted by \ $f:L\xra{\simeq} L'$.

In \cite[II,\S 4-5]{QuR}, Quillen defined 
closed model category structures for the categories \ $d\Lc$ \ and \ 
$s\Lie$, \ as well as for topological spaces (and thus for \ $\TTq$),
and proved:
%
%
\begin{prop}\label{pzero}\stepcounter{subsection}
There are pairs of adjoint functors \ $\TTq\rlh s\Lie$ \  and \ 
$s\Lie\rlh d\Lc$, \ which induce equivalences between the
corresponding homotopy categories: \ 
$\ho\TTq\approx\ho(s\Lie)\approx\ho(d\Lc)$. 
\end{prop}

\begin{notation}\label{nlx}\stepcounter{subsection}
To every simply-connected space \ $\X\in\TTo$ \ one can thus 
associate a DGL \ $\lrp{\Ls,\partial_{L}}\in d\Lc$, \ unique up to
quasi-isomorphism, which determines its rational homotopy type.  
We denote any such DGL by \ $\LX$. \ In particular, \ 
$\His\lrp{\LX}\cong \pim{\X}\otimes_{\ZZ}\Q$, \ the rational homotopy 
algebra of $\X$, which we denote by \ $\Lsx\in\LL$.
\end{notation}

\begin{defn}\label{dcof}\stepcounter{subsection}
The graded Lie algebra \ $\His\lrp{\LX}$ \ does not suffice to
determine the rational homotopy type of \ $\X\in\TTo$: \ in fact,
there may be infinitely many DGLs \ $\{L^{(n)}\}_{n=1}^{\infty}$ \
with \ $\His\lrp{L^{(n)}}\cong \His\lrp{\LX}$, \ no two 
of which are quasi-isomorphic as DGLs; see e.g.\ \cite{LSigD}.
We shall denote by \ $d\Lc(\X)$ \ the full subcategory of \ $d\Lc$ \
whose objects $A$ satisfy \ $\His A\cong \His\lrp{L_{X}}$, \ with the 
isomorphism in $\LL$ \ (see \cite{SStasD}, \cite{LSigD}, or
\cite{FelDT} for treatments of the cohomology analogue of \ $d\Lc(\X)$ \
in terms of algebraic varieties). \ 

The objects of \ $\ho\,d\Lc(\X)$ \
are thus all rational homotopy types which are indistinguishable from \ 
$\XQ$ \ on the primary homotopy operation level. 
Among these there is a distinguished simplest one:
recall that a space \ $\XQ\in\TTq$ \ (or its corresponding DGL model \
$\LX\in d\LL$) \ is called \textit{coformal} (cf.\ \cite{MNeisF}) if \
$\LX$ \ is weakly equivalent to the trivial DGL \ $\lrp{\Ls,0}$ \ 
(where of course \ $\Ls=\His\lrp{L_{X}}$).
\end{defn}

\subsection{Minimal models}
\label{smm}\stepcounter{thm}

Baues and Lemaire (in \cite[Cor.\ 2.4]{BLeM}; see also 
\cite[Props.\ 5.6, 8.1 \& 8.8]{NeiL}) showed that
each connected DGL \ $\lrp{\Ls,\partial}$ \ has a \textit{minimal model} \ 
$\lrp{\hLs,\hat{\partial}}$, \ such that \ $\hLs$ \ is a free graded Lie
algebra, \ $\hat{\partial}:\hat{L}\to\hat{L}$ \ factors through \ 
$[\hat{L},\hat{L}]$, \ and there is a quasi-isomorphism of 
DGLs \ $\varphi:\lrp{\hLs,\hat{\partial}} \to \lrp{\Ls,\partial}$ \ 
(unique up to chain homotopy). In particular, we can choose such a
minimal model \ $\hLx$ \ for any space \ $\X\in\TTo$ \ (cf.\ \S
\ref{nlx}).

As Neisendorfer observes in \cite[\S 5]{NeiL}, in general minimal models
do not exist for non-connected DGLs (but see \cite{MeiM} or \cite{GHTanR} 
for ways around this).

\subsection{Bigraded Lie algebras}
\label{sbgl}\stepcounter{thm}
A differential bigraded Lie algebra, or \textit{DBGL}, 
is a bigraded vector space \ 
$\Lss=\oplus_{p=0}^{\infty} \oplus_{s=0}^{\infty} L_{p,s}$, \ 
equipped with a differential \ 
$\partial_{L}=\partial_{L}^{p,s}:L_{p,s}\to L_{p-1,s}$ \
and a bilinear graded product \ 
$[~,~]:L_{p,s}\otimes L_{q,t}\to L_{p+q,s+t}$ \ for each \ $p,q,s,t\geq 0$ \ 
satisfying:

\setcounter{equation}{\value{thm}}\stepcounter{subsection}
\begin{equation}\label{ebgl}
%
\begin{split}
 [x,y]= (-1)& ^{(p+s)(q+t)+1}[y,x]  \\
  (-1)^{(p+s)(r+u)} & [[x,y],z]  + (-1)^{(p+s)(q+t)}[[y,z],x] \\
 &  + (-1)^{(q+t)(r+u)}[[z,x],y] =0 \\ 
  \partial_{L}\circ \partial_{L}  = & 0   \\
   \partial_{L}[x,y] = & [\partial_{L} x,y] + (-1)^{p+s}[x,\partial_{L} y]  
\end{split}
%
\end{equation}
\setcounter{thm}{\value{equation}}

\noindent for \ $x\in L_{p,s}$, \ $y\in L_{q,t}$, \ and \ $z\in L_{r,u}$.\hsm
The category of such DBGLs will be denoted by \ $\db\LL$, \ with \ 
$\db\Lc$ \ the subcategory with \ $L_{p,0}=0$ \ for all $p$.

\begin{defn}\label{dagl}\stepcounter{subsection}
For each DBGL \ $\lrp{\Lss,\partial_{L}}$ \ there is an
\textit{associated} DGL \ $\lrp{\Ls,\partial_{L}}$, \ defined \
$L_{n}=\bigoplus_{p+q=n} L_{p,q}$ \ 
(same \ $\partial_{L}$); \ some authors re-index \ $\Lss$ \ so that \ 
$\hat{L}_{p,s}=L_{p,p+s}$, \ and then \ $\Ls$ \ is obtained from \ 
$\hat{L}_{\ast,\ast}$ \ by disregarding the first (homological) grading.
\end{defn}

As for ordinary graded Lie algebras, one can define closed model
category structures on \ $s\Lc$ \ and \ $\db\Lc$ \ (see 
\cite[\S 2]{BStG}, and \cite[\S 4]{BlaN}), \ and we have the following
analogue of \cite[I, Props.\ 2.3 \& 4.6, Thm.\ 4.4]{QuR}:
%
%
\begin{prop}\label{pone}\stepcounter{subsection}
There are adjoint functors \ $s\Lc\oulra{N}{N^{\ast}} \db\Lc$, \ which
induce equivalences of the corresponding homotopy categories \ 
$\ho(s\Lc)\approx\ho(\db\Lc)$. \ $N^{\ast}$ \ takes free DBGLs to 
free simplicial graded Lie algebras.
\end{prop}

\begin{proof}
(We give the proof mainly to fix notation which will be needed later.) \ 
Given a simplicial graded Lie algebra \ $\Lds\in s\Lc$, \ we define
the simplicial Lie bracket \ 
$\lrb{~}{~}:L_{p,s}\otimes L_{q,t}\to L_{p+q,s+t}$ \ on \ $\Lds$ \ 
by combining the Lie brackets with the simplicial structure
on \ $\Lds$ \ via the Eilenberg-Zilber map: \ 

\setcounter{equation}{\value{thm}}\stepcounter{subsection}
\begin{equation}\label{eone}
\lrb{x}{y} = 
 \sum _{(\sigma,\tau)\in S_{p,q}} \ 
(-1)^{\varepsilon(\sigma)+pt} 
[s_{\tau_{q}}\dotsc s_{\tau_{1}} x, s_{\sigma_{p}}\dotsc s_{\sigma_{1}} y]
\end{equation}
\setcounter{thm}{\value{equation}}

\noindent where \ $S_{p,q}$ \ denotes the set of all \
$(p,q)$-shuffles \ -- \ that is, partitions of \
$\{0,1,\dotsc,p+q-1\}$ \ into disjoint sets \ 
$\sigma=\{\sigma_{1},\sigma_{2},\dotsc,\sigma_{p}\}$ \ and \ 
$\tau=\{\tau_{1},\tau_{2},\dotsc,\tau_{q}\}$ \ with \ 
$\sigma_{1}<\sigma_{2}<\dotsb <\sigma_{p}$, \ 
$\tau_{1}<\tau_{2}<\dotsb <\tau_{q}$ \ -- \ and \ 
$\varepsilon(\sigma)=p+\sum_{i=1}^{p}(\sigma_{i}-i)$, \ so \ 
$(-1)^{\varepsilon(\sigma)}$ \ is the sign of the permutation
corresponding to \ $(\sigma,\tau)$. \ (See \cite[VIII, \S 8]{MacH}).

Now let \ $\lrp{\Css,\partial}$ \ be the Moore chain complex (cf.\ 
\cite[\S 22]{MayS}) of \ $\Lds$, \ defined by:
\setcounter{equation}{\value{thm}}\stepcounter{subsection}
\begin{equation}\label{ezero}
C_{p,s}=\bigcap_{i=1}^{p}[\Ker(d_{i}^{p})]_{s}\hsp \text{with} \ \ 
\partial_{p}=(-1)^{s}\,d_{0}^{p}\rest{C_{p,s}}.
\end{equation}
\setcounter{thm}{\value{equation}}

%
%
\begin{lemma}\label{lone}\stepcounter{subsection}
If \ $\Ad\in s\LL$ \ is a simplicial graded Lie algebra, \ $x\in A_{p}$ \  
with \ $d_{i}x=0$ \ for \ $1\leq i\leq p-1$, \ and \ $y\in A_{q}$ \
with \ $d_{j}y=0$ \ for \ $1\leq j\leq q-1$, \ then \ 
$d_{k}(\lrb{x}{y})=0$ \ for \ $1\leq k\leq p+q-1$.
\end{lemma}

\begin{proof}
By definition \eqref{eone} we have \ 
\setcounter{equation}{\value{thm}}\stepcounter{subsection} 
\begin{equation}\label{eeleven}
\lrb{x}{y}=\sum _{(\sigma,\tau)\in S_{p,q}}
(-1)^{\varepsilon(\sigma)+p\lrv{y}} [s_{\tau_{q}}\dotsc s_{\tau_{1}} x, 
s_{\sigma_{p}}\dotsc s_{\sigma_{1}} y]\in A_{p+q}.
\end{equation}
\setcounter{thm}{\value{equation}}

Now for each summand \ $w_{\sigma,\tau}\DEF [s_{\tau}x,s_{\sigma} y]$ \ 
in \eqref{eeleven}, with \ $(\sigma,\tau)$ \ a \ $(p,q)$-shuffle,
there are two cases to consider:

The first is that there exist \ $\ell,m$ \ such that \ $\tau_{\ell}=k$, \ 
$\sigma_{m}=k-1$ \ -- \ in which case there is an associated \ 
$(p,q)$-shuffle \ $(\sigma',\tau')$, \ differing from 
$(\sigma,\tau)$ \ only in that \ $\tau_{\ell}$ \ and \ 
$\sigma_{m}$ \ are switched, so that 
$d_{k}(w_{\sigma,\tau}) = d_{k}(w_{\sigma',\tau'})$ \ but \ 
$(-1)^{\varepsilon(\sigma)}= -(-1)^{\varepsilon(\sigma')}$, \ and these
pairs thus cancel in the sum \eqref{eeleven}.

In the second case, \ $k,k-1\in\{\sigma_{1},\dotsc,\sigma_{p}\}$, \ say, \ and 
then there is some \ $0\leq \ell\leq q$ \ with \ $\tau_{\ell}<k-1$ \ 
and \ $\tau_{\ell+1}>k$. \ Since necessarily \ 
$k+1-p\leq \ell\leq k-1$,\ we find that \ $d_{k}s_{\tau}x = 
s_{\tau_{q}-1}\dotsb s_{\tau_{\ell+1}-1} s_{\tau_{\ell}}\dotsb s_{\tau_{1}}
d_{k-\ell}x = 0$.
\end{proof}

\begin{cor}\label{cone}\stepcounter{subsection}
$\lrb{~}{~}$ \ restricts to a bracket \ $C_{p,s}\otimes C_{q,t}\to C_{p+q,s+t}$
\end{cor}

Moreover, if we forget the Lie structure, the Moore chain complex functor $N$  
induces an equivalence between the categories of simplicial graded 
vector spaces and bigraded chain complexes (cf.\ \cite[Thm 1.9]{DoH}),
with the inverse functor $\Gamma$ defined for such a chain complex \ 
$\lrp{\Ass,\partial}$ \ by
$$
(\Gamma\Ass)_{n,s}\DEF~\bigoplus_{0\leq \lambda \leq n}~
                 \bigoplus_{I \in \II_{n,\lambda}} A_{n-\lambda,s}
$$
\noindent (where for each \ $n\geq 0$ \ and \ $0\leq\lambda \leq n$, \
we let \ $\II_{\lambda,n}$ \ denote the set of all sequences of 
$\lambda $ non-negative integers \ $i_{1} < \dotsb < i_{\lambda }(<n)$),
with the obvious face maps (induced by $\partial$) and degeneracies
(see  \cite[p.\ 95]{MayS}).

The left adjoint \ $N^{\ast}:\db\Lc\to s\Lc$ \ to $N$ is defined \ 
$$
N^{\ast}(\lrp{\Lss,\partial})= L(\Gamma(\Lss))/I(\Lss),
$$
where $L$ is the free graded Lie algebra functor, \ and \ $I(\Lss)$ \ is the 
ideal generated by \ $\lrb{\Gamma(x)}{\Gamma(y)}-\Gamma([x,y])$. \ 
The identities \eqref{ebgl} follow from the corresponding ones
in the singly-graded case and the simplicial identities.
\end{proof}
%
%
\sect{Simplicial resolutions}
\label{csr}

The proper algebraic setting for defining our higher homotopy
operations is a suitable notion of a simplicial resolution of \ $\piq{\X}$:

\begin{defn}\label{dcuga}\stepcounter{subsection}
Recall that a category of {\em universal graded algebras\/} 
(or variety of graded algebras, in the terminology of \cite[V,\S 6]{MacC}) 
is a category $\C$ in which the objects are graded 
sets \ $\Xs$, \ together with an action of a fixed set of 
$n$-ary graded operators \ $W=\{
\omega:X_{k_{1}}\times X_{k_{2}}\times\dotsb\times X_{k_{n}}\to X_{m}\}$, \ 
satisfying a set of identities $E$, and the morphisms are functions on
the sets which commute with the operators.  Such categories always
come equipped with a ``free graded algebra'' functor \
$F:\gr\Set\to\C$, \ left adjoint to the ``underlying graded set'' functor \ 
$U:\C\to\gr\Set$. \ In all the examples we shall be concerned with,
the objects \ $\Xs$ will be ``underlying-abelian'' (see 
\cite[\S 2.1.1]{BStG}), and in fact will have the underlying 
structure of a graded vector space over $\Q$.

Examples include $\LL$, and the categories of associative 
(resp.\ non-associative) graded algebras. Note that any ordinary
ungraded category of universal algebras may be thought of as a CUGA
with all objects concentrated in degree $0$.  
\end{defn}

\begin{defn}\label{dsr}\stepcounter{subsection}
A \textit{free simplicial resolution} of an object $B$ in a CUGA \ $\C$
is a weak equivalence from a cofibrant object \ $\Ad\in s\C$ \ to the
constant simplical object associated to $B$ (with respect to the 
closed model category structure on the category \ $s\C$ \ defined in
\cite[II, \S 4]{QuH}). 
Such resolutions always exist, by \cite[II, \S 4]{QuH}; \ see section 
\ref{cmr} below for a specific construction.
\end{defn}

\subsection{Bisimplicial objects}
\label{sbo}\stepcounter{thm}

We shall be interested in a particular type of simplicial resolution,
which may be defined for an arbitrary CUGA \ $\C$ \ ((cf.\
\cite{DKStE} and \cite{BStG}), though we shall
only need it for the case where $\C$ is a category of ungraded
universal algebras, such as \ $\Lie$ \ or \ $\Alg$:

Consider the category \ $ss\C$ \ of bisimplicial objects over $\C$. \ 
We think of an object  \ $\Add\in ss\C$ \ as having 
\emph{internal} and \emph{external} simplicial structures, 
with corresponding homotopy group objects \ $\pii{t}\Add$ \ and \ 
$\pie{s}\Add$ \ (each taking value in \ $s\C$ \ -- \ see \cite[App.]{BStG}). \ 
Let \ $sF:\gr\Ss\to s\C$ \ denote the free graded algebra functor,
extended dimensionwise, and let \ $S^{n}(k)_{\bullet}$ \ be the graded
simplicial set having the simplicial $n$-sphere \ 
$S^{n}_{\bullet}\DEF\Delta[n]/\Delta[n]^{n-1}$ \ in degree $k$. \ 
We think of the simplicial graded algebras \ $F(S^{n}(k)_{\bullet})$ \ as the
$\C$-spheres, or \textit{models}, for \ $s\C$ \ (cf.\ \cite[\S 3.1]{BStG}).
(In the ungraded case one can of course omit the extra degree $k$, and 
write simply \ $F(\bS{n})$). \ Similarly, if \ $D^{n}(k)_{\bullet}$ \
is the graded simplicial set having \ $\Delta[n]$ \ in degree $k$, we
can think of \ $F(D^{n}(k)_{\bullet})$ \ as the $\C$-disks for \
$s\C$. \ The full subcategory of \ $s\C$ \ whose objects are weakly equivalent
to coproducts of such models will be denoted by \ $\M_{\C}$, \ or simply $\M$.

One can use these models to define the so-called  \ ``$E^{2}$-model 
category structure'' for \ $ss\C$, \  as in \cite[\S 5]{DKStE}, \ 
in which  a map \ $f:\Xdd\to \Ydd$ \ is a weak equivalence if 
\setcounter{equation}{\value{thm}}\stepcounter{subsection}
\begin{equation}\label{etwo}
f_{\star}:\pi_{s}\pii{t}\Xdd\to \pi_{s}\pii{t}\Ydd
\text{ \ is an isomorphism for each \ }s,t\geq 0
\end{equation}
\setcounter{thm}{\value{equation}}

We shall not need an explicit description of the fibrations and
cofibrations in \ $ss\C$, \ but only a particular type of cofibrant
object, as follows:

\begin{defn}\label{dmfr}\stepcounter{subsection}
A bisimplicial object \ $\Add\in ss\C$ \ is called \textit{$\M$-free}
if for each \ $m\geq 0$ \ there are graded simplicial sets \ 
$X[m]_{\bullet}\simeq \bigvee _{i}\bS{n_{i}}(k_{i})_{\bullet}$ \ such that \ 
$A_{\bullet,m}\cong F(X[m]_{\bullet})$\ (so that \ $A_{\bullet,m}\in\M$), \ 
and the external degeneracies of \ $\Add$ \ are induced under $F$ by maps \ 
$X[m]_{\bullet}\to X[m+1]_{\bullet}$ \ which are, up to homotopy, 
the inclusion of sub-coproduct summands. \  Any \ $\Xd\in s\C$ \ 
may be resolved by an $\M$-free bisimplicial algebra \ $\Add$ \ 
(see \cite[\S 4.1]{BStG}); this is called an 
\textit{$\M$-free resolution} of \ $\Xd$.
\end{defn}

\begin{defn}\label{ddiag}\stepcounter{subsection}
The \textit{diagonal} of a bisimplicial object \ $\Add\in ss\C$ \ 
is a simplical object \ $\diag(\Add)\in s\C$ \ with \ 
$\diag(\Add)_{n}\DEF A_{n,n}$, \ face maps \ 
$d_{k}=d_{k}^{i}\circ d_{k}^{e}$, \  and degeneracies \ 
$s_{j}=s_{j}^{i}\circ s_{j}^{e}$.
\end{defn}

\begin{remark}\label{rbso}\stepcounter{subsection}
There is a first quadrant spectral sequence with \ 
$$
E^{2}_{s,t}=\pi_{s}^{e}(\pi_{t}^{i}\Add) \Rightarrow
\pi_{s+t}\diag(\Add)
$$
(see \cite{QuS}, and compare \cite[Thm B.5]{BFrH}).

Thus in particular if \ $\Add\to \Xd$ \ is a resolution (in the \ 
$E^{2}$-model category sense), we see that \ 
$\varepsilon :A_{0,\bullet}\to \Xd$ \ induces a  weak equivalence \ 
$\diag(\Add)\simeq \Xd$.

Moreover, the same is true if we disregard the degeneracies and consider 
only the $\Delta$-bisimplicial
resolution \ $\Add^{\Delta}\to \Xd$.
\end{remark}
%
%
\sect{Resolutions for rational spaces}
\label{crrs}

Given a simply-connected space \ $\X\in\TTo$, \ the first 
approximation to an algebraic description of its rational homotopy type is 
given by its rational homotopy Lie algebra \ $\Lsx\DEF\pim\XQ\in\LL$. \ 

If \ $\XQ$ \ were coformal \ (\S \ref{dcof}), then in particular all 
higher homotopy operations vanish in \ $\pis\XQ$, \ and no information beyond \
$\Lsx$ \ itself is needed to determine the rational homotopy type of $\X$. \ 
The higher homotopy operations we shall describe may thus be thought of as 
``obstructions to coformality'', much in the spirit (though not 
the specific approach) of \cite{HStaO}.
 
\subsection{Topological resolutions}
\label{str}\stepcounter{thm}

To proceed further, we need some kind of a ``topological'' simplicial
object \ $\Cd$ \ which realizes a suitable ``algebraic'' simplicial 
resolution \ $\Vds\to\Lsx$ \ in \ $s\LL$, \ in the sense that \ 
$\Vds=\pim\Cd$. \ 
The higher homotopy operations we want then arise as the obstructions 
to realizing the ``algebraic'' augementation map \ $\pim\Cd\to\Lsx$ \ 
topologically.

This can be done using actual topological spaces, as in the integral case 
(see \cite[\S 7]{BlaHH}, as simplified in \cite[\S 4.9]{BlaHO}), but for 
rational spaces it is more convenient to use an algebraic model, in a 
category such as \ $d\LL$. \ To allow us freedom in choosing this model, 
we give a general definition:

\begin{assume}\label{ahho}\stepcounter{subsection}
Let \ $\gr\C$ \ be a CUGA (which we may assume to 
have the underlying structure of a graded vector space), \ and $\C$
the category of (ungraded) universal algebras corresponding to objects
of \ $\gr\C$ \ concentrated in degree $0$. 
The cases we shall be interested in are \ $\C=\Lie$ \ (with \
$\gr\C=\LL$) \ and \  \ $\C=\Alg$ \ (with \ $\gr\C=\A$).

As shown in \cite[App.]{BStG}, for each simplicial algebra \ $\Ad\in s\C$, \ 
the graded homotopy object \ $\pis\Ad$ \ actually takes value in \ $\gr\C$.

For a given  \ $\Ad\in s\C$, \ let \ $\Cdd\to \Ad$ \ be an 
$\M_{\C}$-free resolution (Definition \ref{dmfr}). In particular, this
implies that upon applying the functor \ $\pis$ \ we obtain a free 
simplicial resolution \ $\pii{\ast}\Cdd$ \ (in the ``external''
direction!) of the graded algebra \ $\pis \Ad$. \ In fact, we only
need a $\Delta$-bisimplicial resolution (\S \ref{rbso}), but we shall
nevertheless usually abuse notation by writing \ $\Cdd$ \ for \ $\CdD$.

Next, assume we are given another object \ $\Bd\in s\C$, \ together with 
an isomorphism \ $\varphi:\pis \Ad\cong\pis \Bd$ \ (in \ $\gr\C$). \ 
Define a sequence of morphisms \ $\psi_{n}:\pis \Cnd\to\pis \Bd$ \ by \ 
$\psi_{0}\DEF\varphi\circ\varepsilon_{\#}$ \ and \ 
$\psi_{n+1}\DEF\psi_{n}\circ d_{0}$ \ (which implies that \ 
$\psi_{n+1}=\psi_{n}\circ d_{i}$ \ for all \ $0\leq i\leq n$, \ by the
simplicial identities).

We choose once and for all a fixed map \ 
$f_{0}:C_{0,\bullet}\to\Bd$ \ realizing \ $\psi_{0}$ \ (this is
possible because \ $\Cdd\to\Ad$ \ is $\M$-free) \ 
and define \ $f_{n}:\Cnd\to \Bd$ \ inductively by setting \ 
$f_{n+1}\DEF f_{n}\circ d_{n}$, \ so that \ $\pis(f_{n})=\psi_{n}$ \ 
for all \ $n\geq 0$. \ It is usually most convenient to set \
$f\rest{\DD^{k}}=0$ \ for all $\C$-disks \ $\DD^{k}\hra
C_{0,\bullet}$ \ (cf.\ \S \ref{sbo}).

Note that, because \ $\Cdd$ \ is $\M$-free, the maps \
$\{\psi_{n}\}_{n=0}^{\infty}$ \ define an augmented $\Delta$-simplicial 
object \ $^{h}\CdD\to\Bd$ \ \textit{in the homotopy category} \ 
$\ho(s\C)$ \ -- \ or equivalently, an augmented $\Delta$-simplicial object 
up-to-homotopy.
\end{assume}

\begin{defn}\label{dns}\stepcounter{subsection}
Let \ $D[n]\in\Sa$ \ denote the standard simplicial $n$-simplex, \ 
together with an indexing of its non-degenerate $k$-dimensio\-nal faces \ 
$D[k]^{(\gamma)}$ \ by the composite face maps \ 
$$
\gamma=d_{i_{n-k+1}}\circ\ldots\circ d_{i_{n}}:\mathbf{n}\to\mathbf{k}
$$  
in \ $\Do$ \ (cf.\ \cite[\S 4]{BlaHO}). \ Its \ $(n-1)$-skeleton, which is a
simplicial \ $(n-1)$-sphere, is denoted by \ $\partial D[n]$. \ 
We shall take \ $\ast\DEF D[0]^{(d_{0}d_{1}d_{2}\dotsc d_{n-1})}$ \ as the
base point of \ $D[n]\in\Sa$, \ and we choose once and for all a 
fixed isomorphism \ $\varphi^{(\gamma)}:D[k]^{(\gamma)}\to D[k]$ \ for 
each face \ $D[k]^{(\gamma)}$ \ of \ $D[n]$ \ (see, e.g., \cite[(4.5)]{BlaHO}).
\end{defn}


\begin{defn}\label{dhs}\stepcounter{subsection}
Given \ $\Yd\in s\C$ \ and a simplicial set \ $\Kd\in\Ss$, \ 
we define their \textit{half-smash} \ (in \ $s\C$) \ by: \ 
$$
\Yd\rtimes \Kd\DEF\Yd\otimes \Kd/(\{0\}\otimes\Kd)
$$
(where \ $(\Yd\otimes \Kd)_{n}\DEF \coprod_{x\in K_{n}}~(Y_{n})_{(x)}$ \ -- \
cf.\  \cite[II, \S 1, Prop.\ 2]{QuH}).

Similarly, the \textit{smash product} (in \ $s\C$) \ of \ $\Yd$ \ with a 
\textit{pointed} simplicial set \ $\Kd\in\Sa$ \ is defined \ 
$\Yd\wedge \Kd\DEF\Yd\rtimes\Kd/(\Yd\rtimes\{\ast\})$, \ and if \ 
$\Kd=\bS{r}$ \ (the simplicial sphere), \ we write \ $\Sigma^{r}\Yd$ \
for \ $\Yd\wedge\bS{r}$. \ 
\end{defn}

\begin{remark}\label{rsusp}\stepcounter{subsection}
If \ $\Yd=F(\bS{n})$ \ is a $\C$-sphere (see \S \ref{sbo}),
then \ $\Sigma^{r}\Yd\cong F(\bS{n+r})$ \ is also a $\C$-sphere. \
In fact, many of the usual properties of spheres in \ $\ho\TT$ \ also
hold for $\C$-spheres \ -- \ e.g., \ $\pi_{r}\Xd\cong
[F(\bS{n}),\Xd]_{s\C}$ \ for any \ $\Xd\in s\C$ \ (cf.\ 
\cite[I, \S 4]{QuH}), \ and \ 
$\Yd\simeq \coprod_{i}F(\bS{n_{i}})\Rightarrow 
\Sigma^{r}\Yd\simeq \coprod_{i}F(\bS{n_{i}+r})$ \
(cf.\ \cite[I, \S 3]{QuH}).
\end{remark}

\begin{defn}\label{dcs}\stepcounter{subsection}
Under the assumptions of \S \ref{ahho}, for each \ $n\in\N$, \ we define a \ 
$\partial D[n]$-\textit{compatible sequence} to be a sequence of maps \ 
$\{ h_{k}:C_{k,\bullet}\rtimes D[k]\to \Bd \}_{k=0}^{n-1}$, \ 
such that \ $h_{0}=f_{0}$ \ (under the natural identification \ 
$C_{0,\bullet}\rtimes D[0]=C_{0,\bullet}$), \ and for any iterated
face maps \ $\delta=d_{i_{j+1}}\circ\dotsb\circ d_{i_{n}}$ \ and \ 
$\gamma=d_{i_{j}}\circ\delta$ \ ($0\leq j<n$) \ we have
%
\setcounter{equation}{\value{thm}}\stepcounter{subsection}
\begin{equation}\label{eseven}
h_{j}\circ (d_{i_{j}}\rtimes id) = h_{j+1}\circ
(id\rtimes\iota^{\gamma}_{\delta}) \ \ \ \text{on} \
C_{j+1,\bullet}\rtimes D[j], 
\end{equation}
\setcounter{thm}{\value{equation}}
where \ $\iota^{\gamma}_{\delta}:D[j]\to D[j+1]$ \ is the composite \ 
$\iota^{\gamma}_{\delta}\DEF
\varphi^{\delta}\circ\iota\circ(\varphi^{\gamma})^{-1}$. \ Here \ 
$\varphi^{\gamma}$ \ and \ $\varphi^{\delta}$ \ are the isomorphisms
of Definition \ref{dns}, \ and \ 
$\iota:D[j]^{(\gamma)}\to D[j+1]^{(\delta)}$ \ is the inclusion \ 
(compare \cite[Def.\ 4.10]{BlaHO}). 

A sequence of maps \ 
$\{ h_{k}:C_{k,\bullet}\rtimes D[k]\to \Bd \}_{k=0}^{\infty}$ \
satisfying condition \eqref{eseven} for all $\gamma$, $\delta$, and $n$ 
is called a $\partial D[\infty]$-\textit{compatible sequence}. 
\end{defn}

\begin{defn}\label{dim}\stepcounter{subsection}
Given such a $\partial D[n]$-compatible sequence \ 
$\{ h_{k}:C_{k,\bullet}\rtimes D[k]\to \Bd \}_{k=0}^{n-1}$ \
the \textit{induced map} \ $\bar{h}:\Cnd\rtimes\partial D[n]\to\Bd$ \ 
is defined on the ``faces'' \ $\Cnd\rtimes D[n-1]^{(d_{i})}$ \ of \ 
$\Cnd\rtimes D[n]$ \ by: \ 
$\bar{h}\rest{\Cnd\rtimes D[n-1]^{(d_{i})}}=h_{n-1}\circ (d_{i}\rtimes id)$. \ 
The compatibility condition \eqref{eseven} above guarantees that 
$\bar{h}$ is well-defined.
\end{defn}

\begin{defn}\label{dhho}\stepcounter{subsection}
For each \ $n\geq 2$, \ the $n$-{\em th order homotopy operation\/} 
(associated to the choice of \ $\Cdd\to \Ad$ \ in \S \ref{ahho}) \ 
is a subset \ $\llrr{n}$ \ of the track group \ 
$[\Sigma^{n-1}\Cnd,\Bd]_{s\C}$ \ defined as follows:

Let \ $T_{n}\subseteq [\Cnd\rtimes\partial D[n],\Bd]_{s\C}$ \ 
be the set of homotopy classes of maps \ 
$\bar{h}:\Cnd\rtimes \partial D[n]\to \Bd$ \ induced as above 
by some \ $\partial D[n]$-compatible collection \ $\{ h_{k}\}_{k=0}^{n-1}$. \ 
Since each \ $\Cnd$ \ is a suspension, up to homotopy, by Remark
\ref{rsusp}, we have a splitting 
%
\setcounter{equation}{\value{thm}}\stepcounter{subsection} 
\begin{equation}\label{eeight}
\Cnd\rtimes\partial D[n]\simeq 
(\bS{n-1}\wedge \Cnd) \amalg \Cnd
\end{equation}
\setcounter{thm}{\value{equation}}

\noindent (as for topological spaces). We define \ 
$\llrr{n}\subseteq [\Sigma^{n-1}\Cnd,\Bd]_{s\C}$ \ 
to be the image under the resulting projection of the subset \ 
$T_{n}\subseteq [\Cnd\rtimes\partial D[n],\Bd]_{s\C}$.

Note that the projection of a class \ $[\bar{h}]\in T_{n}$ \ on the other
summand \ $[\Cnd, \Bd]_{s\C}$ \ coming from the splitting 
(\ref{eeight}) is just the homotopy class of the map \ $f_{n}$ \ of 
\S \ref{ahho}. \ 
On the other hand, since \ $\Cdd$ \ was assumed to be $\M$-free, each \
$\Cnd\simeq
\coprod _{k=1}^{\infty }\coprod _{{x\in K_{n,k}}}F(\bS{k}_{(x)})$ \ 
is weakly equivalent to a wedge of spheres over some indexing set \ 
$K_{\ast,\ast}$, so \ $\Sigma^{n-1}\Cnd\simeq
\coprod_{k=1}^{\infty}\coprod_{{x\in K_{n,k}}}F(\bS{k+n-1}_{(x)})$. \ 
Thus \ 
%
\setcounter{equation}{\value{thm}}\stepcounter{subsection} 
\begin{equation}\label{enine}
[\Sigma^{n-1}\Cnd,\Bd]_{s\C}\cong
\prod_{k=1}^{\infty}~\prod_{{x\in K_{n,k}}}~[F(\bS{k+n-1}_{(x)}),\Bd]_{s\C},
\end{equation}
\setcounter{thm}{\value{equation}}
and we shall denote the components of \ $\llrr{n}$ \ under this
product decomposition by \ 
$\llrr{n,x}\subseteq [F(\bS{k+n-1}_{(x)}),\Bd]_{s\C}=\pi_{k+n-1}\Bd$.
\end{defn}

\begin{defn}\label{dcv}\stepcounter{subsection}
An operation \ $\llrr{k}$ \ \textit{vanishes} if it contains the null
class. We say that all the lower order operations \ $\llrr{k}$ \ ($2\leq k<n$)
\textit{vanish coherently} (cf.\ \cite[Def.\ 5.7]{BlaHH}) \ if the \
$\partial D[m]$-compatible collections \  $\{ h_{k}^{\gamma}\}_{k=0}^{m-1}$ \ 
for the various faces $\gamma$ of \ $\partial D[n]$ \ can be chosen to agree on
their intersections, so that they in fact fit together to form a 
$\partial D[n+1]$-compatible collection  \ $\{h_{k}\}_{k=0}^{n}$. \ 
\end{defn}

%
%
\begin{prop}\label{pnine}\stepcounter{subsection}
A necessary condition for the subset \ $\llrr{n}$ \ to be non-empty
that the lower order operations \ $\llrr{k}$ \ vanish for \ $2\leq k<n$; \ 
a sufficient condition is that they vanish coherently. 
\end{prop}

\begin{proof}
See \cite[Theorem 3.29]{BMarkH}.
\end{proof}

\begin{remark}\label{rhho}\stepcounter{subsection}
The coherent vanishing of all the operations \ 
$\{\llrr{n}\}_{n=2}^{\infty}$ \ is equivalent, by 
\cite[Cor.\ 4.21 \& Thm.\ 4.49]{BVoHI} and \cite[\S 4.11]{BlaHH}, to the
\textit{rectifiability} of the augmented $\Delta$-simplicial object
up-to-homotopy \ $^{h}\CdD\to\Bd$: \ that is, its replacement
by augmented $\Delta$-simplicial object \ $\hCd\to\Bd$ \ 
over \ $s\C$ \ (with the simplicial identities now holding precisely,
in \ $s\C$, \ rather than just in \ $\ho(s\C)$), \ such that \ 
$\Cnd\simeq \hat{C}_{n,\bullet}$ \ for each $n$.

This in turn implies (by \S \ref{rbso}) that \ $\diag(\hCd)\simeq\Bd$; \ 
but since \ $\diag(\hCd)\simeq\diag(\CdD)$, \ and \ 
$\diag(\CdD)\simeq \Ad$ \ by assumption, we conclude that \
$\Ad\simeq\Bd$ \ if and only if the higher homotopy operations \ 
$\{\llrr{n}\}_{n=2}^{\infty}$ \ vanish coherently.
\end{remark}

\begin{summary}\label{shho}\stepcounter{subsection}
This yields a first approximation to Theorem A, which may be described as
follows:

We work in \ $\C=\Lie$ \ (and \ $\gr\C=\LL$). \ Given a space \
$\X\in\TTo$ \ we consider the simplicial Lie algebra \ $\Bd$ \
corresponding to a DGL model \ $\LX\in d\LL$ \ for \ $\XQ$ \ (under the
functors of Proposition \ref{pzero}), \ and let \
$\Lsx\DEF\pim\XQ\in\LL$ \ be its rational homotopy Lie algebra, \ 
with \ $\Ad\in s\Lie$ \ the simplicial Lie algebra corresponding to
the trivial DGL \ $L^{(0)}\DEF\lrp{\Lsx,0}$. \ Choose some 
$\M_{\Lie}$-free resolution \ $\Cdd\in ss\Lie$ \ of \ $\Ad$.

$\X$ is coformal if and only if \ $\Ad\simeq\Bd$, \ and this happens if
and only if all the higher homotopy operations \
$\{\llrr{n}\}_{n=2}^{\infty}$ \ associated to \ $\Cdd$ \ vanish
coherently, by Remark \ref{rhho}. \ If not, let \ $n_{0}$ \ denote 
the least \ $n\geq 2$ \ such that \ $0\not\in\llrr{n}$. \ 

Note that we can apply the above procedure to any DGL in \ 
$d\LL(\X)$ \ (Def.\ \ref{dcof}), \ not only to \ $\LX$; \ and the 
existence and vanishing or non-vanishing of the higher homotopy operation \ 
$\llrr{n_{0}}\subset\pis\XQ$ \ is a homotopy invariant. \ 
Denote by \ $\HH{1}$ \ the set of all
homotopy types in \ $\ho\,d\Lc(\X)$ \ for which \ $\llrr{n_{0}}$ \ is
defined and has the same value as for \ $\Bd$ \ itself (i.e., those
DGLs which are indistinguishable from \ $\LX$ \ as far as the primary homotopy
operations, and all the higher homotopy operations \
$\{\llrr{n}\}_{n=2}^{n_{0}}$ \ associated to \ $\Cdd$, \ can see).
For each \ $\alpha\in \HH{1}$, \ choose a representative \ DGL \ 
$L^{(1,\alpha)}$. 

Next, choose a new $\M$-free resolution for the simplicial Lie algebra 
corresponding to \ $L^{(1,\alpha)}$, \ and repeat the above procedure,
yielding a set of higher homotopy operations \ 
$\llrr{n_{1,\alpha}}\subset\pis\XQ$ \ which serve as obstructions to
the existence of a homotopy equivalence \ 
$L^{(1,\alpha)} \xrightarrow{\simeq} \LX$. \ For each such higher 
operation \ $\llrr{n_{1,\alpha}}$, \  we denote by \ $\HH{2,\alpha}$ \ 
the set of all homotopy types in \ $\HH{1}\subseteq ho\,d\Lc(\X)$ \ 
for which \ $\llrr{n_{1,\alpha}}$ \ has the same value as for \ 
$\LX$. \ Now choose representatives \ $L^{(2,\alpha,\alpha')}$ \
for each \ $\alpha'\in \HH{2,\alpha}$, \ and proceed as above. 

In this way we obtain a tree \ $T_{X}$ \ of rational homotopy types in \ 
$\ho d\Lc(\X)$, \ which also indexes a collection of higher homotopy 
operations of the form \ 
$\llrr{n_{k,\alpha_{1},\alpha_{2},\dotsc,\alpha_{k}}}\subseteq\pis\XQ$, \ 
and \ $\lim_{k\to\infty} n_{k}=\infty$ \ along any
branch of the tree \ $T_{X}$, \ so that in fact this collection of operations 
determines the rational homotopy type of $\X$.
\end{summary}

In \cite{BlaAI}, we show how this tree of homotopy types in \
$\ho d\Lc(\X)$, and thus the corresponding collection of higher
homotopy operations, may be described more effectively in terms of a
``Postnikov tower'' for an \ $\M_{\Lie}$-free resolution for $\X$.

%
%
\sect{Minimal resolutions}
\label{cmr}

We now explain how the bisimplicial theory described in section
\ref{crrs} translates into a differential graded theory, when \ 
$\C=\LL$. \ In particular, this allows an application of the
Halperin-Stasheff-F\'{e}lix perturbation theory to our context\vs.

First, it is sometimes convenient to have \textit{minimal} $\M$-free 
resolutions for a DGL, defined for any CUGA \ $\C$ as follows:

\begin{defn}\label{dcwr}\stepcounter{subsection}
Any \ $B\in\C$ \ has a special kind of free simplicial resolution
(see \S \ref{dsr}) \ $\Ad\to B$, \ called a $CW$-\textit{resolution},
defined as follows \ (cf.\ \cite[\S 5.3]{BlaD}): \ 

Assume that for each \ $n\geq 0$, \ $A_{n}=F(\Ts^{n})$ \ is the free
graded algebra on the graded set \ $\Ts^{n}$, \ and that the
degeneracies of \ $\Ad$ \ take \ $\Ts^{n}$ \ to \ $\Ts^{n+1}$. \ 
Let \ $\bar{A}_{n}$ \ denote the sub-algebra of \ $A_{n}$ \ 
generated by the non-degenerate elements in \ $\Ts^{n}$. \ Then we
require that \ $d_{i}\lvert_{\bar{A}_{n}}=0$ \ for \ $1\leq i\leq n$. \
The sequence \ $\bar{A}_{0}=A_{0},\bar{A}_{1},\dotsc,\bar{A}_{n},\dotsc$ \ 
is called a $CW$-\textit{basis} for \ $\Ad$, \ and \ 
$\bar{d}_{0}=d_{0}\lvert_{\bar{A}_{n}}$ \ is the \textit{attaching map} 
for \ $\bar{A}_{n}$.

Such a \ $\Ad\to B$ \ will be called \textit{minimal} if each \ 
$\bar{A}_{n+1}$ \ is minimal among those free algebras in $\C$
which map onto the Moore $n$-cycles \ $Z_{n}\Ad=\Ker(\partial_{n})$ \ 
(see \eqref{ezero}).
\end{defn}

\begin{defn}\label{dls}\stepcounter{subsection}
When \ $\C=\LL$, \ the category of graded Lie algebras, it will be 
more convenient at times to use of the adjoint functors of
Proposition \ref{pone} to replace \ $\Add\to\Xd$ \ by 
a simplicial DGL \ $\Lds\to\Xs$. \ In this case the simplicial models 
are replaced by the corresponding DGLs, namely

\begin{enumerate}
\item A \textit{$d\LL$-$n$-sphere}, \ $\Sph{n}{x}$, \ is 
a DGL of the form \ $\lrp{\bL{\Xs},0}$ \ where \ $\Xs$ \ is the 
graded set with \ $X_{n}=\{x\}$ \ and \ $X_{i}=\emptyset$ \ for \ $i\neq n$.

\item A \textit{$d\LL$-$(n+1)$-disk}, denoted \ $\Di{n+1}{x}$, \ is the DGL \
$\lrp{\bL{\Xs},\partial_{L}}$ \ where \ $X_{n+1}=\{x\}$, \ 
$X_{n}=\{\partial_{L}x\}$, \ and \ $X_{i}=\emptyset$ \ for \ 
$i\neq n,n+1$. \ Its \textit{boundary} is the $d\LL$-$n$-sphere \
$\partial\Di{n+1}{x}\DEF\Sph{n}{\partial_{L}x}$. 

\item A \textit{two-stage DGL} is a DGL \ 
$\lrp{\bL{\Xs},\partial_{L}}\in d\LL$, \ 
where for some \ $n\geq 0$ \ we have \ $X_{i}=\emptyset$ \ for \
$i\neq n,n+1$. \ Any coproduct (in \ $d\LL$) \ of two-stage DGLs will
be called a \textit{free} DGL. 
\end{enumerate}
\end{defn}

Evidently $d\LL$-spheres and disks are free DGLs, and any free DGL may be
described as the coproduct of $d\LL$-spheres and disks \ -- \ more
precisely, as a coproduct of $d\LL$-spheres, disks, and collections of
disks with their boundaries identified to a single sphere. 

\begin{defn}\label{dflr}\stepcounter{subsection}
Following Stover, we define a comonad \ $F:d\LL\to d\LL$ \ by setting
%
\setcounter{equation}{\value{thm}}\stepcounter{subsection}
\begin{equation}\label{ethree}
F(B)\DEF(\coprod_{k=1}^{\infty}~~\coprod_{x\in B_{k}}~~\Di{k}{x})/\sim,
\end{equation}
\setcounter{thm}{\value{equation}}

\noindent for any \ $B=\lrp{\Bs,\partial_{B}}\in d\LL$, \ where we set
\ $\Di{k}{x}\DEF \Sph{k}{x}$ \ if \ $\partial_{B}x=0$, \ and let \ 
$\partial\Di{k+1}{x}\sim\Sph{k}{\partial_{B}x}$ \ if \
$\partial_{B}x\neq 0$.

Clearly \ $F(B)$ is a free DGL, and by iterating $F$ we obtain a free
simplicial DGL \ $\Wds$ \ with \ $W_{n}=F^{n+1}(B)$ \ 
(see \cite[App., \S 3]{GodeT}), 
which we call the \textit{canonical free simplical DGL resolution} of \ 
$B=\lrp{\Bs,\partial_{B}}$, \ which we denote by \ $\Wds(B)$. \ 
Observe that \ $\Wds$ \ (or equivalently, the corresponding
bisimplicial Lie algebra \ $\Wdd$) \ is an $\M$-free resolution of $B$. 
\end{defn}

\begin{remark}\label{rflr}\stepcounter{subsection}
Note that if \ $\partial_{B}\equiv 0$, \ by definition \eqref{ethree} \ 
$F(\Bs,0)$ \ has only spheres, and no disks, and thus the canonical
resolution \ $\Wds(B)$ \ has \ $\partial_{W_{n}}=0$ \ for all \ $n\geq 0$. \ 
Thus \ $\Wds$ \ may be identified with the usual canonical resolution
of the graded Lie algebra \ $\Bs$ \ (coming from the ``free graded Lie
algebra on underlying graded set'' comonad), which we shall denote by \ 
$\Vds(\Bs)$. \ 

However \ -- \ unlike the canonical resolution \ -- \  the construction 
above can be mimicked topologically (cf.\ \cite[\S 2.3]{StoV}).  Since
we want to present our results in a manner which could be generalized
(as far as possible) to the integral case, we have chosen the somewhat
convoluted description of \eqref{ethree}.

Note further that by \eqref{etwo}, if we apply the functor \
$\His$ \ to \ $\Wds\to \Bs$ \ -- \ or equivalently, the functor \ 
$\pii{\ast}$ \ to \ $\Wdd\to\Bd$ \ -- \ we obtain a free simplicial
resolution of the graded Lie algebra \ $\Ls\DEF\His\lrp{\Bs,\partial_{B}}$.
\end{remark}

\begin{notation}\label{ncr}\stepcounter{subsection}
If we write \ $\lrc{x}\in F(B)$ \ for the generator corresponding to an
element \ $x\in \Bs$, \ then recursively a typical DGL generator for \
$W_{n}=W_{n,\ast}$ \ (in the canonical resolution \ $\Wds(B)$) \ is \
$\lrc{\alpha}$, \ for \ $\alpha\in W_{n-1}$, \ so an
element of \ $W_{n}$ \ is a sum of iterated Lie products of elements 
of \ $\Bs$, \ arranged within \ $n+1$ \ nested pairs of brackets \
$\lrc{\lrc{\dotsb}}$. \ With this notation, \ the $i$-th face map of \
$\Wds$ \ is ``omit $i$-th pair of brackets'', and the
$j$-th degeneracy map is ``repeat $j$-th pair of brackets''. \ The
operation of bracketing is defined to be linear\ $\lrc{-}$ \ is linear \ -- \ 
i.e., we set \  $\lrc{\alpha x+\beta y}=\alpha\lrc{x}+\beta\lrc{y}$ \ for \
$\alpha,\beta\in\Q$ \ and \ $x,y\in B$.
\end{notation}

In order to construct minimal $\M$-simplicial resolutions, 
first consider the coformal case:

\subsection{The bigraded model}
\label{sbm}\stepcounter{thm}

Any coformal DGL (\S \ref{dcof}), and in particular \  $L=\lrp{\Ls,0}$, \ 
has a \textit{bigraded model} \ $\Ass\to\Ls$ \ -- \
that is, a bigraded DGL \ $\lrp{\Ass,\partial_{A}}$ \ (see \S \ref{sbgl}) \ 
which is minimal in the sense of \S \ref{smm} (so in particular free
as a graded Lie algebra), along with a quasi-isomorphism \
$\Ass\to\Ls$. \ The bigraded model is unique up to isomorphism.  See
\cite[Part I]{OukiH} for an explicit construction. 

This is just the Lie algebra version of the bigraded model of 
\cite[\S 3]{HStaO} \ (see also \cite{FelM}), \ which is in turn 
essentially the Tate-Jozefiak resolution (see \cite{JozeT}) of a 
graded commutative algebra. 

$A=\lrp{\As,\partial_{A}}$ \ will denote the DGL associated to \ $\Ass$ \
(Definition \ref{dagl}); \ by construction $A$ is the minimal model
(\S \ref{smm}) for $L$ (which is not minimal itself, unless \ $\Ls$ \
happens to be a \textit{free} graded Lie algebra). 


%
%
\begin{prop}\label{ptwo}\stepcounter{subsection}
Let \ $L=\lrp{\Ls,0}\in d\LL$ \ be a coformal DGL, \ and \ 
$\phi:\Ass\to\Ls$ \ its bigraded model; then there is an \ $\M_{d\LL}$-free
simplicial resolution \ $\Cds\to L$, \ with a bijection \ 
$\theta:\Xss\hra\Cds$ \  between a bigraded set \ $\Xss$ \ of 
generators for \ $\Ass$ \ and the set of non-degenerate (\S \ref{dcwr}) \ 
$d\LL$-spheres  in \ $\Cds$. \ Moreover, \ $\His\lrp{\Cds}$ \ is a
minimal $CW$-resolution of \ $\Ls=\His\lrp{\Ass}$, \ with $CW$ basis
generated by \ $\Image(\theta)$. 
\end{prop}

\begin{proof}
By Proposition \ref{pone} there is a simplicial graded Lie algebra 
resolution \ $\Cds\to\Ls$ \  corresponding to \ $\Ass$, \ and thus 
a weak equivalence of simplicial graded Lie algebras \ 
$\psi:\Cds\to\Vds=\Vds(\Ls)$ \ (see \S \ref{rflr}), \ which is 
one-to-one because \ $\Ass$, \ and thus \ $\Cds$, \ are minimal 
(cf.\ \cite[\S 2]{BLeM}). 

Now let \ $\Wds$ \ be the canonical free simplical DGL resolution of \
$\Ass$; \ the fact that \ $\phi:\Ass\to\Ls$ \ is a quasi-isomorphism
implies that there is a weak equivalence \ $\varphi:\Vds\to\Wds$ \ (as
well as one in the other direction). \ The composite \
$\varphi\circ\psi:\Cds\to\Wds$ \ is again a one-to-one weak
equivalence (by minimality); we may therefore think of \ $\Cds$ \ as a
sub-simplicial object of \ $\Wds$.

Moreover, there is an embedding of bigraded vector spaces \
$\eta:\Ass\to\Cds$ \ (see proof of Proposition \ref{pone}), and thus
another such embedding \ $\theta:\Ass\to\Wds$, \ which may be defined
explicitly as follows (using the notation of \S \ref{ncr}):

For \ $x\in X_{0,\ast}$, \ set \ 
$\theta(x)=\lrc{x}\in C_{0,\ast}=F(\As)$. \ Since $\phi$ maps \ 
$X_{0,\ast}$ \ onto a (minimal) set of Lie algebra generators for \ 
$\Ls=\His\lrp{\Ass}$, \ each \ $\theta(x)$ \ is a \ $\partial_{W}$-cycle, 
so \ $\Cus{0}{0}\DEF
\coprod_{k=1}^{\infty}\coprod_{x\in X_{0,k}}\ \Sph{k}{\theta(x)}$ \ 
is a sub free DGL of \ $W_{0,\ast}$.

By minimality of \ $\Ass$, \ any \ $x\in X_{n,\ast}$ \ ($n\geq 1$) \ is 
uniquely determined by \ $\partial_{A}(x)\in A_{n-1,\ast}$. \ Thus if we 
require $\theta$ to be multiplicative (with respect to the ordinary 
bracket in \ $\Ass$, \ and with respect to the simplicial Lie bracket \ 
$\lrb{~}{~}$ \ of \eqref{eone} in \ $\Wds$), \ we may define \ 
$\theta:\Ass\to\Wds$ \ inductively by 
%
\setcounter{equation}{\value{thm}}\stepcounter{subsection}
\begin{equation}\label{efour}
\theta(x)=\lrc{\theta(\partial_{A}(x))},
\end{equation}
\setcounter{thm}{\value{equation}}
and we shall write \ $x^{(0)}$ \ for \ $\theta(x)$ \ if \ $x\in\Xss$.

By definition (see Proposition \ref{pone}), \ 
$d_{0}\circ\theta =\theta \circ \partial_{A}$, \ so for \ 
$x\in X_{n,\ast}$ \ ($n\geq 2$) \ we have \ 
$d_{1}(x^{(0)}) = d_{1} \lrc{\theta(\partial_{A}(x))} =
\lrc{d_{0}\theta(\partial_{A}(x))} = \lrc{\theta(\partial_{A}^{2}(x))}=0$, \  
while \ $\varepsilon(d_{1}(x^{(0)}))$ \ is a \
$\partial_{A}$-boundary for \ $x\in X_{1,\ast}$ \ 
(where \ $\varepsilon: W_{0,\ast}\to \As$ \ is the augmentation).
Thus Lemma \ref{lone} below implies that for $x\in X_{n,\ast}$ \
($n\geq 1$) \ we have \ 
$d_{i}(x^{(0)})=0$ \ for all \ $1\leq i\leq n-1$, \ while \ 
$d_{n}(x^{(0)})$ \ is a \ $\partial_{W}$-boundary.
 
Therefore, if we set \ $\Cus{0}{n}\DEF
\coprod_{k=1}^{\infty}\coprod_{x\in X_{n,k}}\ \Sph{k}{x^{(0)}}$ \ 
for all \ $n\geq 0$, \ we see \ $\{\His\lrp{\Cus{0}{n}}\}_{n=1}^{\infty}$ \ 
is an $\LL$-CW basis for \ $\His\lrp{\Cds}$\vsm . 

In order to give an explicit description of \ $\Cds$ \ in terms of \ 
$\Cuds{0}$, \ we need to know the Lie disks in which \ $d_{n}(x^{(0)})$ \ 
(and their faces) lie. By a double induction on \ $n\geq 1$ \ and \ 
$1\leq r\leq n$, we shall now define, for all \ $x\in X_{n,k}$, \
elements \ $x^{(r)}\in W_{n-r,k+r}$ \ such that \ 
$\partial_{W}(x^{(r)})=d_{n-r}(x^{(r-1)})$:

Note that for each \ $x\in A_{n,\ast}$ \ we have \ 
$$
\partial_{A}(x)=
\sum_{t} a_{t}\omega_{t}\lro{y_{i_{1}},\dotsc,y_{i_{m_{t}}}},
$$
where \ 
$\omega_{t}\lro{\dotsc}$ \ is some \ $m_{t}$-fold iterated Lie bracket, \ 
$y_{i_{j}}\in X_{n_{j},\ast}$ \ with \ $\sum_{j=1}^{m_{t}}n_{j}=n$, \ 
and \ $a_{t}\in\Q$. \ Then
%
\setcounter{equation}{\value{thm}}\stepcounter{subsection}
\begin{equation}\label{efive}
\theta(x)=  \lrc{\theta(\partial_{A}(x))} = 
\lrc{\sum_{t} a_{t}
\omega_{t}\lrbb{y_{i_{1}}^{(0)},\dotsc,y_{i_{m_{t}}}^{(0)}}},
\end{equation}
\setcounter{thm}{\value{equation}}
\noindent where \ $\omega_{t}\lrbb{\dotsc}$ \ is the same \ 
$m_{t}$-fold iterated Lie bracket as above, but now with respect to 
the simplicial Lie bracket \ $\lrb{\,}{\,}$, \ rather than \ $[\,,\,]$.

If we set \ $x^{(s)}=0$ \ for \ $i>n$, \ we may define \ 
$x^{(s)}$ \ for \ $0<s\leq n$ \ inductively by:
%
\setcounter{equation}{\value{thm}}\stepcounter{subsection}
\begin{equation}\label{etwelve}
x^{(s)}=\lrc{\sum_{t} a_{t} 
\sum_{\substack{r_{1}+\dotsb +r_{m_{t}}=s\\0\leq r_{j}}}
\omega_{t}\lrbb{y_{i_{1}}^{(r_{1})},\dotsc,y_{i_{m_{t}}}^{(r_{m_{t}})}}}
\in \Cu{s}{n-s}{k-n+s}.
\end{equation}
\setcounter{thm}{\value{equation}}

Thus if we assume by induction that we have chosen \ 
$y_{i_{j}}^{(r_{j})}$ \ with \ 
$\partial_{W}(y_{i_{j}}^{(r_{j})})=d_{n_{j_{0}}}(y_{i_{j}}^{(r_{j}-1)})$, \ 
it follows from Lemma \ref{lone} below that \ 
indeed \ $\partial_{W}(x^{(s+1)})=d_{n}(x^{(s)})$ \ and \
$d_{i}(x^{(s)})=0$ \ for \ $0<i<n$.

For example, \ $y^{(0)}=\lrc{y}$ \ and \
$\varepsilon(\lrc{y})=y\in A_{k}$ \ for any \ $y\in X_{0,k}$. \ 
Therefore, for \ $x\in X_{1,\ast}$ \ we have \ 
$\varepsilon d_{1}(x^{(0)})=\varepsilon d_{0}(x^{(0)})=\partial_{A}(x)$, \ 
so we may set \ $x^{(1)}=\lrc{x}\in W_{1,k+1}$, \ with \ 
$\partial_{W}(x^{(1)})=d_{1}(x^{(0)})$.

Now if we define by induction
$$
\Cus{r}{n}\DEF\Cus{r-1}{n}\ \amalg \coprod_{k=1}^{\infty} 
\coprod_{x\in X_{n+r,k}}\ \Di{k+r}{x^{(r)}},
$$
\noindent then it is not hard to see that \ $\Cds$ \ is 
the sub-simplicial graded Lie algebra of \ $\Wds$ \ generated (under 
the degeneracies of \ $\Wds$) \ by \ $(\Cus{r}{n})_{r=0}^{n}$ \ for
all \ $n\in\N$. \ (In particular, this is closed under face maps and 
includes \ $\Image(\theta)$, \ and \ $\theta:\Ass\to\Cds$ \ is a weak
equivalence. \ The only non-degenerate Lie spheres in \ $\Cds$ \ are
those of \ $\Cuds{0}$, \ as required.
\end{proof}


\subsection{The filtered model}
\label{sfm}\stepcounter{thm}

If \ $B=\lrp{\Bs,\partial_{B}}\in d\LL$ \ is an arbitrary DGL, it no
longer has a bigraded model, as in \S \ref{sbm} above.  However, it
does have a \textit{filtered model}, constructed as follows:

Let \ $\Ls\DEF\His\lrp{B}$ \ be the homotopy Lie algebra of \ $\Bs$, \ 
and \ $\lrp{\Ass,\partial_{A}}$ \ the  bigraded model for \
$\lrp{\Ls,0}$. \ The filtered model for $B$ is the free graded Lie
algebra \ $\Ass$, \ equipped with an increasing filtration \ 
$0=\F{-1}{A}\subset \F{0}{A}\subset\dotsb \F{r}{A}\subset 
\F{r+1}{A}\subset \dotsb$ \ (defined by \ 
$\F{r}{A}\DEF\bigoplus_{i=0}^{r}A_{i,\ast}$), \ and a new differential \ 
$D_{A}=\partial_{A}+\delta_{A}$ \ on \ $\Ass$ \ such that \ 
$\delta_{A}:A_{n,\ast}\to \F{n-2}{A}$ \ (Of course, \ $D_{A}$ \ is
still required to be a derivation).

We may decompose \ $D_{A}:A_{n,\ast}\to \Ass$ \ 
as \ $D_{A}=\partial_{0}+\partial_{1}+\dotsb+\partial_{n-1}$, \ where \ 
$\partial_{r}:A_{n,\ast}\to A_{n-r-1,\ast}$ \ (and \
$\partial_{0}=\partial_{A}$, \ the original differential of the
bigraded model).

See \cite[II]{OukiH} or \cite{HaralP}; this is again the Lie algebra 
version of a construction of Halperin-Stasheff and F\'{e}lix
in \cite[\S 4]{HStaO},\cite{FelDT}.  

Note that the filtered model is no longer unique, since
its construction depends on choices; in particular, it is not
necessarily minimal.  One again has the associated DGL \
$\lrp{\As,D_{A}}$, \ which is quasi-isomorphic to the original
DGL $B$, \ and \ $\Ass$ \ is obtained by filtering \ $\As$.
%
%
\begin{prop}\label{pthree}\stepcounter{subsection}
Let \ $B=\lrp{\Bs,\partial_{B}}$ \ be a DGL, \ and \ $\lrp{\Ass,D_{A}}$ \ 
a filtered model for $B$; then there is an \ 
$\M_{d\LL}$-free  simplicial resolution \ $\Eds\to B$, \ with a bijection \ 
$\theta:\Xss\hra\Eds$ \  between a bigraded set \ $\Xss$ \ of 
generators for \ $\Ass$ \ and the set of non-degenerate \ $d\LL$-spheres
in \ $\Eds$.
\end{prop}

\begin{proof}
We start with the minimal \ $\M$-free resolution \ $\Cds\to\Ls$ \ for \ 
$\Ls=\His\lrp{\Bs}$, \ constructed as in the proof of Proposition \ref{ptwo},
and deform it into an $\M$-free resolution for
$B$, using the filtered model \ $\lrp{\Ass,D_{A}}$ \ as a guideline.
This time we shall embed the resulting $\M$-free resolution in
the canonical free DGL resolution \ $\Wds$ \ of \ $\lrp{\As,D_{A}}$, \
the DGL associated to the filtered model:

For each \ $x\in X_{n,k}$ \ (where \ $\Xss$ \ is a bigraded set of
generators for the bigraded Lie algebra \ $\Ass$, \ as above),  set \ 
$x^{(n)}=\lrc{x}\in W_{0,k}$, \ and let \ 
$D_{A}(x)=\partial_{0}(x)+\partial_{1}(x)+\dotsb+\partial_{n-1}(x)$ \
as above, with
$$
\partial_{r}(x)= 
\sum_{t} a_{t}^{(r)}\omega_{t}^{(r)}\lro{y_{i_{1}},\dotsc,y_{i_{m_{t}}}}
\in A_{n-r-1,\ast},
$$
where \ $\omega_{t}^{(r)}\lro{\dotsc}$ \ is some \ 
$m_{t}$-fold iterated Lie bracket, as above, and each \ 
$y_{i_{j}}\in X_{n_{j},\ast}$ \ with \ $\sum_{j=1}^{m_{t}}n_{j}=n-r-1$.

If we set \ $x^{(s)}=0$ \ for \ $i>n$, \ we may define \ 
$x^{(s)}$ \ for \ $0<s\leq n$ \ inductively by:
%
\setcounter{equation}{\value{thm}}\stepcounter{subsection}
\begin{equation}\label{esix}
x^{(s)}=\lrc{\sum_{r=0}^{s}\sum_{t} a_{t}^{(r)} 
\sum_{\substack{r_{1}+\dotsb +r_{m_{t}}=s-r\\0\leq r_{j}}}
\omega_{t}^{(r)}\lrbb{y_{i_{1}}^{(r_{1})},\dotsc,y_{i_{m_{t}}}^{(r_{m_{t}})}}}
\in \Cu{s}{n-s}{k-n+s}
\end{equation}
\setcounter{thm}{\value{equation}}

Using Lemma \ref{lone} and the fact that for any \ 
$\Ad\in s\LL$, \ $x\in A_{p}$ \ and \ $y\in A_{q}$ \ we have \ 
$d_{p+q}(\lrb{x}{y})=\lrb{d_{p}(x)}{y} + (-1)^{q}\lrb{x}{d_{q}(y)}$, \
one may then verify inductively that \ 
$d_{n-s}(x^{(s)})=\partial_{W}(x^{(s+1)})$ \
and \ $d_{i}(x^{(s)})=0$ \ for \ $0<i<n-s$, \ for \ all \ $0\leq s<n$.
The rest of the construction is as in the proof of Proposition \ref{ptwo}.
\end{proof}


We have the following analogue of Definition \ref{dhs}:

\begin{defn}\label{ddhs}\stepcounter{subsection}
Given a DGL \ $L=\lrp{\bL{\Xs},\partial_{L}}\in d\LL$ \ and a 
simplicial set \ $A\in\Ss$, \ we define their \textit{half-smash} \ 
$L\rtimes A=\lrp{\bL{\Ys},\partial'}\in d\LL$ \ by setting \ 
$Y_{n}\DEF \coprod_{k=0}^{n}X_{k}\times \hat{A}_{n-k}$, \ 
where \ $\hat{A}_{i}$ \ denotes the set of non-degenerate $i$-simplices
of $A$. \ For \ $a\in\hat{A}_{k}$ \ and \ $x\in X_{m}$, \ we set 
$$
\partial'(x,a)=\sum_{i=0}^{k}(-1)^{i+m} (x,d_{i}a) + (\partial_{L}x,a)
$$
(and extend \ $\partial'$ \ by requiring that it be a derivation).
\end{defn}


\begin{remark}\label{rot}\stepcounter{subsection}
In order to apply the obstruction theory described in \S \ref{shho}, 
note that all the definitions of section \ref{crrs} pass over to
the DGL setting in a straightforward manner. However, if we now start
with the trivial DGL \ $A=\Ls^{(0)}\DEF\lrp{\Lsx,0}$, \ we may take \ 
$\Cds\to \Ls^{(0)}$ \ to be the minimal $\M$-free resolution of
Proposition \ref{ptwo}, corresponding to the bigraded model for \ 
$\lrp{\Lsx,0}$, \ and let \ $B=\lrp{\Bs,\partial_{B}}$ \
(corresponding to \ $\Bd$ \ in \S \ref{shho}) be the filtered model 
for \ $\LX$. \ We assume that \ $A\not\simeq B$.

As explained in \S \ref{shho}, there is a least \ $n_{0}\geq 2$ \ 
such that \ $0\not\in\llrr{n_{0}}\subseteq \His\lrp{B}$, \ and we write \ 
$\llrr{n_{0}}=(\llrr{n_{0},x_{i}})_{i\in I}$, \ in the notation of
\ref{dhho}, where \ $x_{i}\in X_{n_{0},i}$ \ and \ 
$\Sph{k_{i}}{x_{i}}$ \ are corresponding DGL spheres in \ $C_{n_{0},\ast}$ \ 
(we include in the index set $I$ only those coordinates of \eqref{enine}
which do not vanish).

Again let \ $\HH{1}$ \ denote the set of all homotopy types in \ 
$\ho d\LL(\X)$ \ for which \ $\llrr{n_{0}}$ \ has
the same value as for \ $\LX$, \ and choose a representative \ 
$L^{(1,\alpha)}\in d\LL(\X)$ \ for each \ $\alpha\in \HH{1}$. \ By 
\cite[\S 3]{HStaO}, we may assume \ $L^{(1,\alpha)}$ \ is obtained
from $B$ by perturbation of \ $\partial_{B}$. \ Proceeding as in \S
\ref{shho} we obtain a tree of DGLs \ 
$L^{(k,\alpha_{1},\alpha_{2},\dotsc,\alpha_{k})}\in d\LL(\X)$, \ and 
by \cite[Theorem 3.1]{BlaD}, we know that \ 
$L^{(k,\alpha_{1},\alpha_{2},\dotsc,\alpha_{k})}$ \ may be chosen 
to agree with \ $\LX$ through degree \ $n+1$ \ at least, so \ 
$\colim_{n} L^{(n)}\simeq\LX$ \ along any branch of the tree.

Note also that because  \ $\His\lrp{\Cds}$ \ is a (minimal) CW resolution 
of \ $\His\lrp{B}$, \ in each case, the maps \
$\psi_{n}:\His\lrp{C_{n,\ast}}\to \His\lrp{B}$ \ are null for all \ 
$n\geq 1$ \ (see \S \ref{ahho}). \ Thus any \ 
$\partial D[n]$-compatible collection \ $\{ h_{k}\}_{k=0}^{n-1}$ \ in
\S \ref{dhho} induce a map \ $C_{n,\ast}\wedge\partial D[n]\to B$ \
directly, without need of the splitting \eqref{eeight}.
\end{remark}


\begin{remark}\label{rrep}\stepcounter{subsection}
The second order operation described in the previous example is
actually a secondary Whitehead product. Unlocalized higher order Whitehead
products were defined by G.\ Porter in \cite[1.3]{GPorW},
and the relation between this definition and the rational version has
been studied by several authors \ -- \ see \cite{AArkS}, 
\cite{AlldR,AlldR2}, \cite{RetaL,RetaMO} and \cite[V.1]{TanrH}.
However, there are other higher order rational homotopy operations, too: \ 

For example, in the DGL \ 
$\Ls=\lrp{\bL{a_{1},b_{1},c_{1},d_{1},x_{4},y_{4},z_{4},w_{4}},\partial}$, \ 
with \ $\partial(x)=[[b,a],c]$, \ $\partial(y)=[[b,a],d]$, \ 
$\partial(z)=[[d,c],a]$ \ and \ $\partial(w)=[[d,c],b]$, \ the cycle \
$[x,d]+[y,c]+[z,b]+[w,a]$ \ represents such an operation. There
appears to be no general procedure for representing these as 
\textit{integral} higher order operations in \ $\pis \X$; \ we shall
offer a (partial) answer to this difficulty in section \ref{cna}.
\end{remark}

%
%
\sect{Homology of DGLs}
\label{chd}

Obstructions in algebraic topology traditionally take values in
suitable cohomology groups. In order to show that this holds in our
setting, too, we recall Quillen's definition of homology and cohomology
in model categories:

\begin{defn}\label{dao}\stepcounter{subsection}
An object $X$ in a category $\C$ is said to be \textit{abelian} if it is an 
abelian group object  \ -- \ that is, if \ $\Hom_{\C}(Y,X)$ \ has a natural 
abelian group structure for any \ $Y\in \C$. \ 
When $\C$ is \ $\Lie$, \ $\Alg$, \ $s\Lie$, \ $s\Alg$, \ $\LL$, or \ 
$d\LL$, \  for example, this is equivalent to requiring that all products 
vanish in $X$ \ (cf.\ \cite[\S 5.1.3]{BStG}). 

The full subcategory of abelian objects
in $\C$ is denoted by \ $\C_{ab} \subset \C$. \ In the cases of
interest to us, this will itself be an abelian category. It is equivalent
to the category \ $\Vect$ \ of vector spaces if \ $\C=\Lie$ \ or \ $\Alg$, \ 
to $\VV$ if \ $\C=\LL$, \ to the category \ $s\Vect$ \ of simplicial 
vector spaces if \ $\C=s\Lie$ \ or \ $s\Alg$, \ and to the category \
$d\VV$ \ if \ $\C=d\LL$ \ (see \S \ref{snac}).

In these cases,  we have an \textit{abelianization} functor \ 
$Ab:\C\to \C_{ab}$, \ along with a natural transformation \ 
$\theta:Id\to Ab$ \ having the appropriate universal property.
In all the examples above, \ 
$Ab(X)=X/I(X)$, \ where \ $I(X)$ \ is the ideal in \ $X\in\C$ \
generated by all non-trivial products. 
\end{defn}

\begin{defn}\label{dhom}\stepcounter{subsection}
Let $\C$ be a category as above, which also has a closed model category
structure: \ in \cite[II, \S 5]{QuH} (or \cite[\S 2]{QuC}), 
Quillen defines the \textit{homology} of an object \ 
$X\in\C$ \ to be the total left derived functor \ $\LLL(Ab)$ \ of \ 
$Ab$, \ applied to $X$ \ (cf.\ \cite[I, \S 4]{QuH}).

In more familiar terms, this means that we construct a resolution \
$A\to X$ \ (i.e., replace $X$ by a weakly equivalent cofibrant object \ 
$A\in\C$), and then define the $i$-th homology group of $X$ by \ 
$H_{i} X\DEF \pi_{i}(Ab(A))$, \ for an appropriate concept of 
homotopy groups \ $\pis$ \ in \ $\C_{ab}$ \ (see \cite[II, \S 4]{QuH}). \ 
One must verify, of course, that this definition is independent of the
choice of the resolution \ $A\to X$.

Similarly, the \textit{cohomology} of $X$ with coefficients in \ 
$M\in\C_{ab}$ \ is defined:
$$ \ 
H^{i}(X;M)\DEF [\LLL(Ab)X,\Omega^{i+N}\Sigma^{N}M]_{\C_{ab}} \hsp 
\text{for \ }N\hsp \text{large enough}
$$
(where the loop and suspension functors $\Omega$ and $\Sigma$ are 
defined in  \cite[I, \S 2]{QuH}).

Again, in the cases that interest us, $\Omega$ is essentially
the shift operator \ $\Sigma^{-1}$ \ of \S \ref{snac}, and so
the $i$-th cohomology group of $X$ with coefficients in $M$ is then \ 
$H^{i}(X;M)\DEF [\Sigma^{i}Ab(A),M]_{\C_{ab}}$. 
\end{defn}

\begin{defn}\label{dhl}\stepcounter{subsection}
If $\C$ itself does not have a closed model category structure, one
often defines the homology of\ $\X\in\C$ \ by embedding $\C$ in some
category which does have such a structure, which in many cases may be
taken to be \ $s\C$, \ the category of simplicial objects over $\C$ \ 
(see \cite[II, \S 4]{QuH}). Thus, if \ $\iota:\C\hra s\C$ \ is 
the embedding of categories defined by taking \ $\iota(C)$ \ to be 
the constant simplicial object equal  to $C$ in all dimensions, \ then \ 
$H_{i}(C)\DEF \pi_{i}(\LLL(Ab\circ \iota) C)$.

This is the approach usually taken for \ $\C=\Lie$, \ $\Alg$, \ 
or $\LL$: \ to define the homology of a graded Lie algebra \ 
$\Ls\in\LL$, \ say, one chooses a free simplicial resolution \ 
$\Ads\to\Ls$ \ (such as the canonical resolution \ -- \ cf.\ 
\S \ref{dflr}), and then calculate the homotopy groups of the
simplicial graded vector space \ $Ab(\Ads)\in s\VV$ \ (or the homology
groups of the bigraded chain complex in \ $\db\VV$ \ corresponding to \
$Ab(\Ads)$ \ -- \ see proof of Proposition \ref{pone}).
\end{defn}

\begin{remark}\label{rhl}\stepcounter{subsection}
Note that if we apply Definition \ref{dhom} as is to a DGL \ 
$L=\lrp{\Ls,\partial_{L}}\in d\LL$, \ we may take the resolution 
$A$ to be the minimal model \ $\hat{L}=\lrp{\hLs,\hat{\partial}}$ \ 
for $L$ (cf.\ \S \ref{smm}), \ and since its abelianization is just 
the graded vector space \ $Q(\hat{L})$ \ of indecomposables, and \ 
$Q(\hat{\partial})=0$ \ by definition, \ \ $H_{i}(L)$ \ would be 
isomorphic to the vector space spanned by a set of generators for 
$\hat{L}$ in dimension $i$. 

If we want cohomology with coefficients in an object \ 
$\Ms\in d\LL_{ab}\approx d\VV$ \ with \textit{trivial} differential \ -- \ 
i.e., \ $\Ms$ \ is just a graded vector space \ -- \ we find \ 
$$
H^{i}(X;\Ms)\cong \prod_{j=1}^{\infty}\ \Hom_{\Vect}(H_{j}(X),\ M_{i+j}), 
$$
by the universal coefficients theorem.
\end{remark}

However, since $L$ is itself graded, we would like \ $H_{\ast}L$ \ 
to be bigraded (with a ``homological'' degree, as well as a
``topological'' one). This requires a combination of Definitions 
\ref{dhom} and \ref{dhl}, as follows:

\begin{defn}\label{dhdgl}\stepcounter{subsection}
The \textit{homology} \ $\Hss(\Ld)$ \ of a simplicial Lie algebra \ 
$\Ld\in s\Lie$ \ is defined to be the left 
derived functors of the abelianization, with respect to the 
$E^{2}$-closed model category structures (\S \ref{sbo}) on \ $ss\Lie$ \ and \ 
$ss\Lie_{ab}\approx dd\VV$ \ respectively. More precisely, 

\setcounter{equation}{\value{thm}}\stepcounter{subsection} 
\begin{equation}\label{eten}
H_{s,t}(\Ld)\DEF \pi_{s}(\LLL(Ab\circ \iota)\Ld)_{t}
=\pi_{s}\pii{t}(Ab\Add),
\end{equation}
\setcounter{thm}{\value{equation}}
where \ $\Add\to\Ld$ \ is some $\M$-free bisimplicial resolution of \
$\Ld$.

Similarly, for any DGL \ $L\in d\LL$ \ we may define \ 
$H_{s,t}\lrp{L}\DEF \pi_{s}\Hi{t}\lrp{Ab(\Ads)}$, \ for a \ $\M_{d\LL}$-free 
simplicial resolution \ $\Ads\to L$; \ and these two
definitions agree under the equivalence of homotopy categories \ 
$\ho(s\Lie)\approx\ho(d\LL)$ \ of Proposition \ref{pzero}.

The bigraded cohomology of a DGL $L$ with coefficients in the abelian
DGL (i.e., chain complex) $M$ is defined analogously as \
$$
H^{s}_{t}\lrp{L;M}\DEF \pi^{s}(\Hom_{d\LL_{ab}}(Ab(\Ads), M)_{t}).
$$
\end{defn}

We note that the homology and cohomology of differential graded
(commutative) algebras have been defined by Goodwillie 
(in \cite{GooC}) and Burghelea \& Vigu\'{e}-Poirrier (in
\cite{BVigC}), in a manner analogous to the traditional definitions of
Hochschild homology. See  \cite[\S 5.3]{LodC}.

%
%
\begin{prop}\label{peight}\stepcounter{subsection}
For any DGL \ $L\in d\LL$, \ there is a monomorphism of graded
vector spaces \ $H_{n,t}\lrp{L}\hra H_{n,t}\lrp{L'}$, \ where \
$L'\simeq \lrp{\His\lrp{L},0}$ \ is the coformal model for $L$; the
same holds for cohomology with trivial coefficients.
\end{prop}

\begin{proof}
If \ $\Ass$ \ is the bigraded model for \ $L'$, \ and \ $\Cds\to L'$ \
the simplicial resolution of Proposition \ref{ptwo}, then 
the non-degenerate spheres \ $\Sph{k}{x^{(0)}}\subset C_{n,t}$, \ 
which correspond to a vector space basis for \ $H_{n,t}\lrp{L'}$, \ 
are in bijective correspondence with the generators \ 
$x\in X_{n,t}$ \ for \ $\Ass$.

Now let \ $\Bss$ \ be a filtered model for $L$ obtained by perturbing \
$\lrp{\Ass,\partial_{A}}$, \ and \ $\Eds\to L$ \ the associated simplicial 
resolution of Proposition \ref{pthree}: \ since \ $\Bss$ \ need no
longer be minimal (\S \ref{sfm}), a vector space basis for \ 
$H_{n,t}\lrp{L}$ \ now corresponds to a \textit{subset} of the collection
of non-degenerate spheres \ $\Sph{k}{x^{(0)}}\subset E_{n,t}$, \ 
(which are still in bijective correspondence with the generators \ 
$x\in X_{n,t}$ \ for \ $\Ass$ \ or \ $\Bss$).
\end{proof}

Note that this description of the homology implies that \ 
$\Hss(L)$ \ is indeed just a bigraded version of the DGL homology 
defined in \S \ref{dhl}.

%
%
\begin{prop}\label{pfour}\stepcounter{subsection}
The collection of higher operations \ 
$\llrr{n_{k,\alpha_{1},\alpha_{2},\dotsc,\alpha_{k}},x}$ \
which determine the rational homotopy type of \ $\X\in\TTo$ \ 
(described in \S \ref{shho} above) are indexed by elements \ 
$x\in H_{n_{k,\alpha_{1},\alpha_{2},\dotsc,\alpha_{k}},t}
(L^{(k,\alpha_{1},\alpha_{2},\dotsc,\alpha_{k})})$ \ 
in the homology of the DGLs of \S \ref{rot}, and take value in the 
cohomology of these DGLs, with \ 
$$
\llrr{n_{k,\alpha_{1},\alpha_{2},\dotsc,\alpha_{k}},x}\subseteq
H^{n_{k,\alpha_{1},\alpha_{2},\dotsc,\alpha_{k}}}_{t+n_{k}-1}
(L^{(k,\alpha_{1},\alpha_{2},\dotsc,\alpha_{k})};\piq{\X}).
$$
\end{prop}

\begin{proof}
We may construct a simplicial resolution \ $\Eds$ \ for each
successive DGL \ $L^{(k)}=L^{(k,\alpha_{1},\alpha_{2},\dotsc,\alpha_{k})}$, \ 
corresponding to the filtered models obtained as perturbations \ 
$\lrp{\Ass,D_{A}}$ \ of the bigraded model \ 
$\lrp{\Ass,\partial_{A}}$ \ for \ $L^{(0)}$, as above. \ The
non-degenerate spheres \ $\Sph{m}{x}=\Sph{m}{x^{(0)}}\subset E_{m,t}$ \ which
index the higher homotopy operations \ $\llrr{m,x}$ \ are thus in 
bijective correspondence with the generators \ $x\in X_{m,t}$ \ for \ $\Ass$.
However, if $x$ is not minimal \ -- \ in the sense that \
$D_{A}(x)\not\in[A,A]$, \ or \ $x+\alpha= D_{a}(y)$ \ for some \ 
$\alpha\in \Ass$ \ and \ $y\in X_{m+1,t}$ \ -- \ then we can construct
a new simplicial resolution \ $\Eds'$ \ of \ $L^{(k)}$ \ in which \ 
$\Sph{m}{x}$ \ has been eliminated (though of course new spheres may
appear in higher simplicial dimensions). By the universal property of
resolutions (i.e., of cofibrant objects in the \ $E^{2}$ model
category for \ $sd\LL$ \ -- \ see \S \ref{sbo}) there is a map of
resolutions \ $\Eds\to\Eds'$, \ and there can be no non-vanishing 
higher operation \ $\llrr{n_{k,\alpha_{1},\dotsc,\alpha_{k}},x}$ \ 
which serves as an obstruction to rectifying the augmentation up-to-homotopy \ 
$\varphi:\Eds\to L^{X}$, \ since \ $\varphi_{n}\rest{\Sph{m}{x}}$ \ 
can be factored through \ $0\in\Eds'\to L^{X}$. \ Thus the only homotopy 
operations which can appear are those corresponding to
\textit{non-trivial} homology classes in \
$H_{\ast}\lrp{L^{(k,\alpha_{1},\alpha_{2},\dotsc,\alpha_{k})}}$. \
These yield the requisite cohomology classes by Universal Coefficients.
\end{proof}

Proposition \ref{peight} thus implies that we may if we like think of
all the higher homotopy operations described in \S \ref{shho} 
(associated to the various deformations of \ $L^{(0)}$) \ as lying in 
one fixed bigraded group \ $H^{\ast}_{\ast}(L^{(0)};\piq{\X})$, \ which
is of course just the usual cohomology of a graded Lie algebra, and 
is easier to compute than the cohomology of a non-trivial DGL.

%
%
\sect{Non-associative algebra models}
\label{cna}

The DGL higher homotopy operations
are unsatisfactory from a topological point of view 
because there is no obvious way to translate them, in general, into 
\textit{unlocalized} topological homotopy operations. We now describe an 
algebraic model for rational homotopy theory which may serve to 
answer this objection.

\subsection{Non-associative graded algebras}
\label{snaa}\stepcounter{thm}
Let $\A$ denote the category of \textit{non-asso\-ci\-at\-ive
gra\-ded al\-geb\-ras}: an object \ $\As\in\A$, \ is just a graded
vector space equipped with a bilinear graded product \ 
$A_{p}\otimes A_{q}\to A_{p+q}$. \ Let \ $\bA{X}$ \ denote the free 
non-associative graded algebra generated by a graded set \ $\Xs$. \ 
As in \S\ref{sdgl}, the functor \ ${\mathbb A}:gr\Set\to\A$ \ factors 
through \ $A:\VV\to\A$.

$d\A$ \ will denote the category of differential graded
non-associative algebras \ $\lrp{\As,\partial_{A}}$, \ called 
\textit{DGNAs}; \ the differential \ $\partial_{A}$ \ must satisfy \ 
$\partial_{A}\circ \partial_{A}=0$ \ and \ 
$\partial_{A}(x\cdot y) = 
\partial_{A}x\cdot y + (-1)^{\lrv{x}}x\cdot \partial_{A}y$, \ as for DGLs.

For simplicity we assume each \ $\As\in \A$, $d\A$ \ is connected \ -- \ 
that is, \ $A_{0}=\{0\}$. \ Again, we have a
category \ $\db\A$ \ of differential bigraded non-associative algebras 
(\textit{DBGNA}s), as in \S \ref{sbgl}. 

As for any CUGA (\S \ref{dcuga}), one can define a closed model 
category structure on \ $s\A$ \ (see \cite[II, \S 4]{QuH}) 
and thus by \cite[\S 4]{BlaN} on \ $\db\A$, \ and one has the expected 
analogues of Propositions \ref{pzero} and \ref{pone}:
%
%
\begin{prop}\label{pfive}\stepcounter{subsection}
There are adjoint functors \ $s\Alg\oulra{N}{N^{\ast}} d\A$, \ which
induce equivalences of the corresponding homotopy categories \ 
$\ho(s\Alg)\approx\ho(d\A)$.
\end{prop}

%
%
\begin{prop}\label{psix}\stepcounter{subsection}
There are adjoint functors \ $s\A\oulra{N}{N^{\ast}}\db\A$, \ which
induce equivalences of the corresponding homotopy categories \ 
$\ho(s\A)\approx\ho(\db\A)$; \ and \ $N^{\ast}$ \ takes free DBGNAs to 
free simplicial non-associa\-tive algebras. 
\end{prop}

\begin{notation}\label{njm}\stepcounter{subsection}
For any \ $\lrp{\Xs,\cdot}\in \A$, \ let \ 
$\npr{x}{y}$ \ denote \ $\tfrac{1}{2}(x\cdot y + \epsi{x}{y} y\cdot x)$. \ 
We then have \ $\npr{y}{x}=\epsi{x}{y}\npr{x}{y}$, \ so \ 
$\lrp{\Xs,\npr{\,}{\,}}$ \ is now a non-associative graded algebra with a 
graded-commutative (or: graded skew-symmetric) multiplication. \
Moreover, any graded derivation $\partial$ on \ $\lrp{\Xs,\cdot}$ \ is
also a derivation with respect to \ $\npr{\,}{\,}$, \ and any morphism
of algebras from \ $\lrp{\Xs,\cdot}$ \ to a graded-commutative algebra
will also respect \ $\npr{\,}{\,}$. \ Therefore we can (and will)
assume that our non-associative algebras are all graded-commutative,
and denote the product by \ $\npr{\,}{\,}$. 

Moreover, if \ $\Ad\in s\A$ \ is a simplicial graded algebra, we 
shall also use the notation \ 
$$
\lrb{x}{y} =  \sum _{(\sigma,\tau)\in S_{p,q}} \ 
(-1)^{\varepsilon(\sigma)+p\lrv{y}} 
[s_{\tau_{q}}\dotsc s_{\tau_{1}} x, s_{\sigma_{p}}\dotsc
s_{\sigma_{1}} y]
$$
for the corresponding simplicial bracket (compare \eqref{eone}).
\end{notation}

\begin{defn}\label{dar}\stepcounter{subsection}
Any simplicial Lie algebra \ $\Ld\in s\Lie$ \ is in particular an 
object in \ $s\Alg$; \ let \ $\iota:d\LL\hra d\A$ \ be the
inclusion functor. Note that even if each \ $L_{n}$ \ is free as a Lie 
algebra, it is not free as a non-associative algebra: a free
simplicial resolution \ $\Jd\to \iota(\Ld)$ \ in the category \ 
$s\Alg$ \ (see \S\ref{dsr}) will be called a \textit{$s\Alg$-model} 
for \ $\Ld$. \ Such models can be constructed functorially, 
for example by a variant of \S  \ref{rflr}.

There is also the analogous concept of a \ \textit{$d\A$-model} \
$\Js\in d\A$ \ of a DGL $L$; \ we can of course translate back and 
forth between these two types of models using Proposition \ref{pfive}.

Since the DGL \ $L=\lrp{\Ls,\partial_{L}}$ \ has an internal grading,
and its $d\A$-model \ $J=\lrp{\Js,\partial_{J}}$ \ is constructed as a
resolution of $L$, it is natural to define a second ``homological''
degree on \ $\Js$, \ so as to have a \textit{filtered} $d\A$-model
(cf.\ \S \ref{sfm}). \ If the DGL is trivial (i.e., \ $\partial_{L}=0$), \ 
a $d\A$-model for \ $\lrp{\Ls,0}$ \  will be a differential bigraded 
non-associative algebra (DBGNA) \ 
$J=\lrp{\Jss,\partial_{J}}$ \ (cf.\ \S \ref{sbm}). 
\end{defn}

\begin{remark}\label{rja}\stepcounter{subsection}
Define a \textit{Jacobi algebra} to be a DGNA \
$J=\lrp{\Js,\partial_{J}}\in d\A$ \ such that \ $\His J\in \LL$. \ In
particular, any $d\A$-model $J$ of a DGL $L$ is a Jacobi algebra,
since \ $\His J\cong\His L$. \ We denote by \ $\JJ\subset d\A$ \ the
full subcategory of Jacobi algebras.
These algebras are clearly related to the strongly homotopy Lie
algebras of \cite{SStasD} \ (see also \cite{LMarkS}), \ though in
general a Jacobi algebra is just a ``Lie algebra up to homotopy''.

Note that even though \ $d\A$ \ itself is a CUGA, \ 
$\JJ$ apparently is not, and it does not inherit many desirable 
properties from \ $d\A$: \ for example, \ $\JJ$ is not closed under the
coproduct in \ $d\A$. \ However, one still has \textit{free} Jacobi
algebras, in the following sense:
\end{remark}
%
%
\begin{lemma}\label{ltwo}\stepcounter{subsection}
There is a functor \ $J:d\VV\to d\A$, \ and a natural transformation \
$\theta:J\to L$ \ such that:
\begin{enumerate}
\renewcommand{\labelenumi}{(\alph{enumi})}
\item \ $J\Vs$ \ is free as an algebra, for any chain complex \ $\Vs$;
\item \ $\theta_{\Vs}$ \ is a surjective quasi-isomorphism;
\item any chain map \ $\varphi:\Vs\to\Ks$ \ (where \ $\Vs\in d\VV$ \
and \ $\Ks\in\JJ$) \ extends to a \ $d\A$ map \ $\hat{\varphi}:J\Vs\to\Ks$.
\end{enumerate}
\end{lemma}

\begin{proof}
As noted above, we may define \ $\Jss=J(\Vs)$ \ by induction on the
homological filtration:

Start with \ $J_{0,\ast}=A(\Vs)$ \ (the free non-associative algebra
on the differential graded vector space \ $\lrp{\Vs,\partial_{V}}$, \ 
with \ $\partial_{J}$ \ extending \ $\partial_{V}$ \ as a derivation),
and let \ $\theta_{0}:J_{0,\ast}\to L(\Vs)$ \ be the obvious
surjection, with \ $K_{0,\ast}=\Ker(\theta_{0})$ \ a two-sided $d\A$-ideal of \
$J_{0,\ast}$.

Choose once and for all some collection of generators \ 
$M_{0}=\{\mu_{i}\}_{i\in I}$ \ for \ $K_{0,\ast}$ \ as a \ 
$J_{0,\ast}$-bimodule: \ for each choice of 
a \ $d\VV$-basis \ $\{x_{\gamma}\}_{\gamma\in\Gamma}$ \ 
for \ $\Vs$ \ -- \ that is, of a graded vector space basis of the form \ 
$\{x_{\alpha},\partial_{V}x_{\alpha},x_{\beta}\}$, \ with \ 
$\partial_{V}x_{\beta}=0$ \ -- \ we may write each \ $\mu_{i}$ \ as some
expression \ $\mu_{i}(x_{\gamma_{i_{1}}},\dotsc,x_{\gamma_{i_{n}}})$. \ 
If we then choose some other \ $d\VV$-basis \ 
$\{x'_{\gamma'}\}_{\gamma'\in\Gamma'}$ \ for \ $\Vs$, \ again we will have \ 
$\mu_{i}'\DEF\mu_{i}(x'_{\gamma'_{i_{1}}},\dotsc,x'_{\gamma'_{i_{n}}})
\in K_{0,\ast}$; \ define \ $J_{1,\ast}$ \ to be the free
non-associative algebra on the DG vector subspace of \ $K_{0,\ast}$ \ 
spanned by all such \ ``canonical operations'' 
$\mu_{i}'$ \ $(i\in I$), \ for all possible choices of $d\VV$-bases \ 
$\{x'_{\gamma'}\}_{\gamma'\in\Gamma'}$ \ for \ $\Vs$. \ 
Again one has the obvious augmentation \ 
$J_{1,\ast}\to J_{0,\ast}$ \ to serve as \ $\partial_{J}$ \ (with \ 
$\partial_{J}\circ\partial_{J}=0$ \ by construction), and one takes 
the kernel \ $K_{1,\ast}$ \ of this augmentation for the next step.

Proceeding in this way we may define the functor $J$ by induction on
the homological filtration; if the collections of operations \ 
$M_{n}$ \ are chosen canonically, the functor itself will be
canonical. Properties (a)-(c) are readily verified.
\end{proof}

\begin{example}\label{ejr}\stepcounter{subsection}
Let \ $L=\lrp{\bL{\Xs},0}$ \ be the trivial free DGL on a graded set \
$\Xs$; \ in this case the canonical DBGNA model \
$\lrp{\Jss,\partial_{J}}=J(\Xs,0)$ \ for $L$ may be described in part 
as follows:

Let \ $J_{0,\ast}=\bA{\Xs}$ \ (the free non-associative algebra on \
$\Xs$). \ Since the Jacobi identity holds in \ $\bL{\Xs}$, \ but not 
in \ $\bA{\Xs}$, \ we have \ 
$\mu(x,y,z)= \npr{x}{\npr{y}{z}}-\npr{\npr{x}{y}}{z}+
(-1)^{qr}\npr{\npr{x}{z}}{y} \in K_{0,\ast}$ \ for all \ 
$x,y,z\in \Xs$. \ Thus  \ $J_{1,\ast}$ \ will be generated as a \ 
$J_{0,\ast}$-bimodule by the image of \ $(J_{0,\ast})^{\otimes 3}$ \
under the \ $\Sigma_{3}$-equivariant multilinear map \ 
$\lambda_{3}:J_{i,p}\otimes J_{j,q}\otimes J_{k,r}\to J_{i+j+k+1,p+q+r}$. \ 
Here the symmetric group \ $\Sigma_{n}$ \ acts on \ 
$J_{\ast,p_{1}}\otimes\dotsb\otimes J_{\ast,p_{n}}$ \ 
by permutations, and on \ $J_{\ast,p_{1}+\dotsb+p_{n}}$ \ via the 
\textit{Koszul sign} homomorphism \ 
$\varepsilon_{I}:\Sigma_{n}\to\{1,-1\}$ \  (defined by 
letting \ $\varepsilon_{I}((k,k+1))=(-1)^{i_{k} i_{k+1}+1}$ \ for any adjacent 
transposition \ $(k,k+1)\in \Sigma_{n}$). \ We set
$$
\partial_{J}(\lamt{x}{y}{z})=
\npr{x}{\npr{y}{z}} - \npr{\npr{x}{y}}{z} + 
(-1)^{qr}\ \npr{\npr{x}{z}}{y}.
$$

However, there are relations among these elements \ $\lamt{x}{y}{z}$, \ 
so we define a $\Sigma_{4}$-equivariant multilinear map \ 
$\lambda_{4}:J_{i,p}\otimes J_{j,q}\otimes J_{k,r}\otimes J_{\ell,s}\to 
J_{i+j+k+\ell+2,p+q+r+s}$, \ 
\begin{equation*}
\begin{split}
\partial_{J}(\lamf{x}{y}{z}{w}) & =       
        \npr{x}{\lamt{y}{z}{w}} \ +\ \npr{\lamt{x}{y}{z}}{w} \\
      -\ & \eps{z}{w}\ \npr{\lamt{x}{y}{w}}{z} \\
           \ +\ & (-1)^{\lrv{y}(\lrv{z}+\lrv{w})}\npr{\lamt{x}{z}{w}}{y} \\
      -\ & (\lamt{\npr{x}{y}}{z}{w}\ +\ \lamt{x}{y}{\npr{z}{w}}) \\
      +\ & \eps{y}{z}\ (\lamt{\npr{x}{z}}{y}{z}\ \\
      +\ & \lamt{x}{z}{\npr{y}{w}}) \\
      -\ & (-1)^{\lrv{w}(\lrv{y}+\lrv{z})}\ 
           (\lamt{\npr{x}{w}}{y}{z}\ \\
      +\ & \lamt{x}{w}{\npr{y}{z}}).
\end{split}
\end{equation*}

In fact, one can define a sequence of ``higher Jacobi relations'' \
$\lambda_{n}(x_{1}\otimes\dotsb\otimes x_{n})$, \ for all \ 
$n\geq 3$, \ which yield an explicit construction of \ $J(\Xs)$ \ 
for a the free (graded) Lie algebra \ $L{\Xs}$. \ See \cite{ABlaR}, and 
compare \cite[2.1]{LMarkS}.
\end{example}

\subsection{$\A$-homotopy operations}
\label{saho}\stepcounter{thm}

One can now apply the theory of section \ref{crrs} verbatim to any space \
$\X\in \TTo$ \ with \ $\C=\Alg$ \ rather than \ $\Lie$, \ to obtain
a sequence of higher homotopy operations as in \S \ref{shho}
which determining the rational homotopy type of $\X$ \ -- \ the only
difference being that the simplicial resolutions \ $\Cdd$ \ 
of the successive simplicial Lie algebras \ $\Ld^{(k)}$ \ are now 
$\M_{\Alg}$ free  resolutions of \ $\Ld^{(k)}$ \ in \
$ss\Alg$. 

This is the reason that the theory of section \ref{csr} was 
stated for an arbitrary CUGA, rather than specifically for \ $\C=\Lie$. \ 
The reason that our general theory was stated for \textit{simplicial} rather
than \textit{differential graded} universal algebras is that there
seems to be no reasonable version of Proposition \ref{pfive} for an 
arbitrary CUGA.

\subsection{Minimal resolutions}
\label{smr}\stepcounter{thm}

To make the construction more accessible, it is again useful to have 
\textit{minimal} $\M_{\Alg}$ resolutions, as in section \ref{cmr}.
For this purpose, we consider a variant of the above approach:

Even though \ $\JJ$ \ does not inherit a closed model category 
structure from \ $d\A$, \ one may define \textit{models} for $\JJ$, \
in the sense of \S \ref{sbo}, by letting a \textit{Jacobi sphere} be any
$d\A$-model of a $\LL$-sphere (\S \ref{dls}), \ and more generally let \ 
$\M_{\JJ}$ \ denote the full subcategory of $\JJ$ consisting of DGNAs
weakly equivalent to objects in \ $\M_{d\LL}$ \ -- \ i.e., Jacobi
models of DGLs which are (up to homotopy) coproducts of $d\LL$-spheres.

An \ $\M_{\JJ}$ resolution of a DGL \ $L$, which we shall call simply
a \textit{Jacobi resolution},  is then defined to be a free
simplicial resolution of DGNAs \ $\Ads\to \iota(L)$ \ (Def.\
\ref{dsr}), with each \ $A_{n,\ast}\in\M_{\JJ}$. \ Note that such an \
$\Ads\to \iota(L)$ \ is at the same time also an $\M_{\JJ}$-Jacobi 
resolution of the $d\A$-model \ $\Js$ \ of $L$, and it is
usually more convenient to think of it as such.

There is a comonad \ $F:d\A\to d\A$ \ as in 
\eqref{ethree}, which yields the \textit{canonical Jacobi resolution} \ 
$\Uds$ \ for any \ $C\in\JJ$, \ as in \S \ref{dflr}. Again we 
may use the notation of \S \ref{ncr}. 

One also has an analogue of Propositions \ref{ptwo} and \ref{pthree},
as follows:
%
%
\begin{prop}\label{pseven}\stepcounter{subsection}
Let \ $B=\lrp{\Bs,\partial_{B}}\in d\LL$ \ be any DGL, \ and \ 
$\lrp{\Ass,D_{A}}$ \ a filtered model for $B$: \ then there is a 
Jacobi  resolution \ $\Jds\to \iota(B)$, \ with a bijection \ 
$\theta:\Xss\hra\Jds$ \  between a bigraded set \ $\Xss$ \ of 
generators for \ $\Ass$ \ and the set of non-degenerate \ 
$d\A$-spheres in \ $\Jds$.
\end{prop}

\begin{proof}
Let \ $G=G_{B}$ \ be a $d\A$-model for the DGL $B$, \ and \ $\Uds\to G$ \ 
the canonical Jacobi resolution. As in the proof of Proposition \ref{ptwo}, 
we may define a map \ $\theta:\Ass\to\Uds$ \ inductively by the
equation \ $\theta(x)=\lrc{\theta(\partial_{A}(x))}$ \ (compare
\eqref{efour}), and we shall again write \ $x^{(0)}$ \ for \ 
$\theta(x)$ \ if \ $x\in\Xss$ \ (a set of generators for \ $\Ass$), \
and let \ $\Vss$ \ be the bigraded vector space spanned by \
$\theta(\Xss)$. \ For simplicity of notation we consider first the
case where $B$ has trivial differential and \ $\Ass$ \ is bigraded
(with \ $D_{A}=\partial_{A}$).

For each \ $n\in \N$, \ define the sub-DGNA \ $\Jus{0}{n}$ \ of \ 
$U_{n,\ast}$ \ to be \ $J(V_{n,\ast})$, \ in the notation of Lemma \ref{ltwo} \ -- \ 
that is, \ $\Jus{0}{n}$ \ is the coproduct, in \ $\JJ$, \ of a set 
of Jacobi spheres \ $\Sph{k}{x^{(0)}}$, \ one for each generator \ 
$x\in X_{n,k}$ \ of \ $\Ass$. \ By Lemma \ref{ltwo}(c), it is enough 
to define the face and degeneracy maps of \ $\Jds$ \ on each $x$ \ -- \  
where we may use the description of \S \ref{ncr}.

Once again, we want \ $d_{i}(x^{(0)})$ \ to be a \
$\partial_{U}$-boundary for each \ $1\leq i\leq n$; \ but the
analogues of the elements \ $x^{(s)}$ \ of Propositions \ref{ptwo} and
\ref{pthree} are more complicated, so we need some definitions:

For each \ $0\leq s\leq n$, \ let \ $\KK_{n,s}$ \ denote the set of 
all sequences \ $I=(i_{1},\dotsc,i_{s})$ \ of integers \ 
$1\leq i_{1} < \dotsb < i_{s}\leq n$, \ corresponding to the $s$-fold 
face map \ $d_{I}=
d_{i_{1}}\circ\dotsb\circ d_{i_{s}}:\mathbf{n}\to\mathbf{n}-\mathbf{s}$ \ 
in \ $\Do$ \ (compare Definition \ref{dns} and the proof of
Proposition \ref{ptwo}). \ Given \ $I=(i_{1},\dotsc,i_{s})\in\KK_{n,s}$, \ 
for each \ $1\leq j\leq s$ \ let \ 
$I(\hat{j})\DEF (i_{1},\dotsc,\widehat{i_{j}},\dotsc,i_{s})\in\KK_{n,s-1}$ \ 
be obtained from $I$ by omitting the $j$-th entry. \ By repeatedly using
the identity \ $d_{k}d_{m}=d_{m-1}d_{k}$ \ ($k<m$), \ we can find a unique \
$\kappa(j)\in\{1,2,\dotsc,n\}$ \ such that \ 
$d_{\kappa(j)}\circ d_{I(\hat{j})}=d_{I}$. \ 

For each \ $x\in X_{n,k}$, \ \ $0\leq s\leq n$, \ 
and \ $I\in\KK_{n,s}$, \ we want to choose choose an element \ 
$x^{(s;I)}\in \Ju{s}{n-s}{k-n+s}\subset U_{n-s,k-n+s}$ \ by
induction on \ $n-s$, \ starting with \ 
$x^{(0,\emptyset)}\DEF x^{(0)}=\theta(x)$, \ so that:
%
\setcounter{equation}{\value{thm}}\stepcounter{subsection}
\begin{equation}\label{ethirteen}
\partial_{U}(x^{(s;I)})=
\sum_{j=1}^{s}(-1)^{j} d_{\kappa(j)}(x^{(s-1;I(\hat{j}))})
\end{equation}
\setcounter{thm}{\value{equation}}
\noindent for \ $s\geq 1$. \ (The index $s$ is not really needed,
since \ $s=\lrv{I}$, \ but it is useful for keeping the analogy with
the notation of \eqref{etwelve} in mind.)

Note that since \ $d_{0}\circ\theta=\theta\circ \partial_{A}$ \ no
longer holds here (because \ $d_{0}$ \ is a morphism in $d\A$, \ not
in $d\LL$), \ it is not generally true that \ $d_{1}(x^{(0)})=0$ \ for
all \ $x\in\Xss$. \ However, since applying \ $\His$ \
still yields a $CW$ resolution \ $\His\Juds{0}\to \His C=\His B$, \ 
(where $C$ is the $d\A$-model for the DGL \ $B$), \ we know that \ 
$d_{i}(x^{(0)})$ \ must be a $\partial_{U}$-boundary \ for each \
$1\leq i\leq n$. \ Thus we can choose an element \ 
$x^{(1;1)}\in U_{n-1,k-n+1}$ \ with \ 
$\partial_{U}(x^{(1;1)})=d_{1}(x^{(0)})$ \ (a special case of
\eqref{ethirteen}) \ -- \ and in fact \ $x^{(1;1)}$ \ may be expressed
in terms of the ``canonical operations'' \ $\mu_{i}$ \ of Lemma \ref{ltwo}.

Now let \ $\partial_{A}(x)=
\sum_{t} a_{t}\omega_{t}\lro{y_{i_{1}},\dotsc,y_{i_{m_{t}}}}$ \ for \ 
$y_{i_{j}}\in A_{n_{j},\ast}$, \ so \ 
%
\setcounter{equation}{\value{thm}}\stepcounter{subsection}
\begin{equation}\label{efourteen}
\begin{split}
x^{(0)}=& \lrc{\sum_{t} a_{t}\omega_{t}
\lrbb{y_{i_{1}}^{(0)},\dotsc,y_{i_{m_{t}}}^{(0)}}} \\
 =& \lrc{\sum_{t} 
a_{t}\sum_{(J_{1},\dotsc,J_{m_{t}})}(-1)^{\varepsilon_{t}} \omega_{t}
\lro{s_{J_{1}}(y_{i_{1}}^{(0)}),\dotsc,s_{J_{m_{t}}}(y_{i_{m_{t}}}^{(0)})}}
\end{split}
\end{equation}
\setcounter{thm}{\value{equation}}

\noindent as in \eqref{efive},  where the $(n-n_{j})$-multi-index \ 
$J_{j}\subseteq\{0,1,\dotsc,n-1\}$ \ 
is obtained by repeated shuffles, which also determine the sign \
$(-1)^{\varepsilon_{t}}$ \  (see \eqref{eone} ff.). \ Therefore, 
$$
d_{k}(x^{(0)})=
\lrc{\sum_{t}
a_{t}\sum_{(I_{1},\dotsc,I_{m_{t}})}(-1)^{\varepsilon_{t}} 
\omega_{t}\lro{d_{k-1}s_{J_{1}}(y_{i_{1}}^{(0)}),\dotsc,
d_{k-1}s_{J_{m_{t}}}(y_{i_{m_{t}}}^{(0)})}}.
$$

The proof of Lemma \ref{lone} (which is valid in \ $d\A$, \ too) implies
by induction on \ $t\geq 2$ \ that for each summand \ 
$$
v_{t}=\omega_{t}\lro{s_{J_{1}}(y_{i_{1}}^{(0)}),\dotsc,
s_{J_{m_{t}}}(y_{i_{m_{t}}}^{(0)})},
$$
there is exactly one \ 
$1\leq j\leq t$ \ and \ $0\leq \ell\leq k-1$ \ such that \ 
\begin{equation*}
\begin{split}
d_{k-1}(v_{t})= \omega_{t}[&s_{J'_{1}}(y_{i_{1}}^{(0)}),\dotsc,
s_{J'_{j-1}}(y_{i_{j-1}}^{(0)}), \\
& s_{J''_{j}}d_{\ell}(y_{i_{j}}^{(0)}),
s_{J'_{j+1}}(y_{i_{j+1}}^{(0)}),
\dotsc,s_{J'_{m_{t}}}(y_{i_{m_{t}}}^{(0)})],
\end{split}
\end{equation*}
\noindent for suitable multi-indices \ 
$J'_{1},\dotsc,J'_{j-1},J'_{j+1},\dotsc,J'_{m_{t}}$ \ and \ $J''_{j}$.

Since \ $\ell<k$ \ and \ $n_{j}<n$, \ we may assume by induction
that we have defined \ $y_{i_{j}}^{(1;\ell)}\in U_{n_{j}-1,\ast}$ \ 
such that \ $\partial_{U}(y_{i_{j}}^{(1;\ell)})=d_{\ell}(y_{i_{j}}^{(0)})$, \ 
and then let \ $x^{(1;k)}$ \ be
\begin{equation*}
\begin{split}
\langle\sum_{t}\sum_{(J_{1},\dotsc,J_{m_{t}})}(-1)^{\varepsilon{t}}
\omega_{t}[& s_{J'_{1}}(y_{i_{1}}^{(0)}),\dotsc,
s_{J'_{j-1}}(y_{i_{j-1}}^{(0)}),\\
& s_{J''_{j}}(y_{i_{j}}^{(1;\ell)}),
s_{J'_{j+1}}(y_{i_{j+1}}^{(0)}),
\dotsc,s_{J'_{m_{t}}}(y_{i_{m_{t}}}^{(0)})]\rangle.
\end{split}
\end{equation*}

The rest of the construction of the elements \ $x^{(s;I)}$, \ as well
as the generalization to the filtered case, is similar to that in 
the proofs of Propositions \ref{ptwo} and \ref{pthree}.
\end{proof}


We can now summarize the main result of this paper in the following

%
%
\begin{thm}\label{tone}\stepcounter{subsection}
Let $\X$ be a simply connected space, and \ $\Lsx\DEF\pim\XQ\in\Lc$ \ 
its rational homotopy Lie algebra. \ There is a tree \ $T_{X}$ \ of DGLs \ 
$L^{(k,\alpha_{1},\dotsc,\alpha_{k})}$, \ 
starting with \ $L^{(0)}\simeq\lrp{\Lsx,0}$, \ 
and for each branch \ $\alpha_{1},\dotsc,\alpha_{k},\dotsc$ \ of \
$T_{X}$, \ an increasing sequence of positive integers (or $\infty$) \ 
$(n_{k}=n_{k,\alpha_{1},\dotsc,\alpha_{k}})_{k=1}^{\infty}$ \ such that 
\begin{enumerate}
\renewcommand{\labelenumi}{(\alph{enumi})}
\item \ $\His\lrp{L^{(k,\alpha_{1},\dotsc,\alpha_{k})}}\cong\Lsx$;
\item \ The higher homotopy operations \ $\llrr{m}\subset 
H^{m}_{\ast}\lrp{L^{(k,\alpha_{1},\dotsc,\alpha_{k})};\Lsx}$ \
associated to a minimal Jacobi resolution of \ 
$L^{(k,\alpha_{1},\dotsc,\alpha_{k})}$ \ as in \S \ref{dhho}, \ 
vanish for \ $m<n_{k}$.
\item The operation \
$\llrr{n_{k}}=\llrr{n_{k,\alpha_{1},\dotsc,\alpha_{k}}}\subset 
H^{n_{k}}_{\ast}\lrp{L^{(k)};\Lsx}$ \ does not vanish 
(unless \ $n_{k}=\infty$).
\item \ For any \ $\alpha_{k+1}$ \ along the branch of \ 
$(\alpha_{1},\dotsc,\alpha_{k})$, \ the DGL \ 
$L^{(k+1,\alpha_{1},\dotsc,\alpha_{k},\alpha_{k+1})}$ \ may be chosen
to agree with \ $L^{(k,\alpha_{1},\dotsc,\alpha_{k})}$ \ in degrees 
$\leq n_{k}+1$, \ so the sequential colimit \ 
$L^{(\infty)}=\colim_{k}L^{(k)}$, \ 
along any branch, is well defined, and is a DGL model for $\X$.
\end{enumerate}
\end{thm}

The main difference between the construction described here (in \
$s\JJ$) \ and that of Proposition \ref{pthree} is that the ``higher 
order information'' in \ $\Jds$ \ is no longer concentrated in the
last face map \ $d_{n}:J_{n,\ast}\to J_{n-1,\ast}$. \ As a result, the
higher homotopy operations associated to the (minimal) Jacobi resolution
(as in section \ref{crrs}) are true simplicial operations, which can
be translated more directly into topological ones.

\end{document}